\documentclass[12pt]{article}
\usepackage{amssymb,amsmath}

\numberwithin{equation}{section}

\allowdisplaybreaks

\newenvironment{Proof}{\removelastskip\par\medskip
\noindent{\em Proof.} \rm}{\penalty-20\null\hfill$\square$\par\medbreak}

\newenvironment{Proofx}{\removelastskip\par\medskip
\noindent{\em Proof} }{\par}

\usepackage[colorlinks=true, urlcolor=blue,linkcolor=blue, citecolor=blue]{hyperref}

\def\Dom{{\mathord{{\rm Dom}}}}

\def\real{{\mathord{{\rm I\kern-2.8pt R}}}}        
\def\inte{{\mathord{{\rm I\kern-2.8pt N}}}}
\def\PP{{\mathord{{\rm I\kern-2.8pt P}}}}

\def\real{{\mathord{\mathbb R}}}

\def\inte{{\mathord{\mathbb N}}}

\def\var{{\mathrm{{\rm Var \ \!}}}}

\newtheorem{prop}{Proposition}[section]
\newtheorem{lemma}[prop]{Lemma}
\newtheorem{definition}[prop]{Definition}
\newtheorem{corollary}[prop]{Corollary}
\newtheorem{theorem}[prop]{Theorem}
\newtheorem{remark}[prop]{Remark}

\def\Dom{{\mathrm{{\rm Dom }}}}

\def\Ent{{\mathrm{{\rm Ent \ \!}}}}

\DeclareMathOperator{\esssup}{ess\,sup}

\DeclareMathOperator{\essinf}{ess\,inf}

\title{\huge Surface measures and related functional inequalities 
 on configuration spaces} 
\author{Christian Houdr\'e\thanks{Research supported in part by a 
NSF Grant.} 
\and
Nicolas Privault}
\date{December 9, 2003}

\begin{document}

\maketitle 

\begin{abstract} 
 Using finite difference operators, 
 we define a notion of boundary and 
 surface measure for configuration sets under 
 Poisson measures. 
 A Margulis-Russo type identity and a co-area formula 
 are stated with applications to 
 deviation inequalities and functional inequalities, 
 and bounds are obtained on the associated 
 isoperimetric constants. 
\end{abstract} 
 
\normalsize

\vspace{0.5cm}

\noindent {\bf Key words:} Configuration spaces, Poisson measures, 
 surface measures, co-area formulas, isoperimetry. 
\\ 
{\em Mathematics Subject Classification.} 60G57, 60H07, 28A75, 60D05. 
\small

\normalsize

\baselineskip0.7cm

\section{Introduction} 
 Isoperimetry consists in determining sets with minimal 
 surface measure, among sets of given volume 
 measure. 
 In probability theory, isoperimetry is generally formulated 
 by expressing the volume of sets via a probability measure, 
 and surface measures using the expectation of an appropriate 
 gradient norm. 
 Gaussian isoperimetry is a well-known subject, 
 see, e.g., \cite{ledoux} for a review. 
 A notion of surface measure on configuration spaces 
 has been recently introduced in \cite{bpr} using differential operators. 
 Discrete isoperimetry is also possible 
 on graphs and Markov chains, by defining the surface 
 measure of a set $A$ as an average of the number of elements in 
 $A$ that are connected to an element in $A^c$, 
 cf. e.g. \cite{diaconisstroock}, \cite{lawlersokal}, 
 without requiring any smoothness on $A$. 
 In this framework, an isoperimetric result has been 
 obtained in \cite{bobkovgoetze}, Prop.~3.6, 
 for i.i.d. Poisson vectors in $\inte^d$. 
 
 In this paper we consider the problem of isoperimetry 
 on configuration space in finite volume, i.e. on the space 
 $\Omega$ of a.s. finite configurations $\omega = \{ x_1, \ldots , x_n \}$, 
 $n\geq 1$, of a metric space $X$. 
 The configuration space $\Omega$ is equipped with a Poisson measure 
 $\pi$ with intensity $\sigma$, where $\sigma$ is a finite 
 diffuse Borel measure on $X$. 
 Working with the configuration space instead of finite Poisson distributed 
 i.i.d. vectors is similar to working with measurable functions on $\real$ 
 instead of step functions. 
 Each ($\pi$-a.s. finite) configuration 
 $\omega\in \Omega$ has a set of ``forward'' neighbors 
 of the form $\omega \cup \{ x \}$, $x\in \omega^c = X\setminus \omega$, 
 and a set of ``backward'' neighbors 
 of the form $\omega \setminus \{ x \}$, $x\in \omega$. 
 A Markov chain and a graph of unbounded degree 
 can both be constructed on $\Omega$. 
 In the Markov case one adds a point distributed according to 
 the normalized intensity measure 
 to a given configuration. 
 In the graph case, a point chosen at random is removed 
 from a given configuration. 
 Such operations of additions and subtraction of points 
 are also frequently used in statistical mechanics and in connection 
 with logarithmic Sobolev inequalities, see, e.g., \cite{daipra}. 
 Here they allow to construct two notions of neighbor 
 (respectively denoted forward and backward) for a 
 given configuration. 
 It turns out that the graph and Markov kernels 
 are mutually adjoint under the Poisson measure, 
 and we will work with a symmetrized kernel 
 in order to take both the graph and Markov structures 
 into account. 
 We emphasize that it is necessary here to use 
 the graph and Markov approaches simultaneously 
 (i.e. to consider both forward and backward neighbors), 
 since considering only the Markov part or the graph part separately 
 yields trivial values of the isoperimetric constants $h^\pm_p = 0$. 
 In fact the classical discrete isoperimetric results that hold in our setting 
 are those which are valid both in the Markov and graph cases. 
 This notion of neighbors is used to define the inner and outer 
 boundary and the surface measure $\pi_s$ of arbitrary sets of configurations. 
 Isoperimetry and the related isoperimetric constants 
 are then studied by means of co-area formulas. 
 We can define dimension free isoperimetric constants 
$$ 
 h_1 
 = \inf_{0 < \pi (A ) < \frac{1}{2}} 
 \frac{\pi_s (\partial A)}{\pi (A)}, 
$$ 
 and 
$$ 
 h_\infty = \inf_{0 < \pi (A ) < \frac{1}{2}} 
 \frac{\pi (\partial A)}{\pi (A)}. 
$$ 
 Let $\lambda_2 = 1$ denote the 
 optimal constant in the Poincar\'e inequality on configuration space 
 for the finite difference operator $D$. 
 We have $\frac{1}{2} \leq h_1 \leq 2+2\sqrt{\sigma (X)}$, and 
$$ 
 \max \left( 
 \frac{1}{\sqrt{\pi \sigma (X)}} , 
 \frac{1}{2 \sigma (X)} \right) 
 \leq h_\infty \leq 
 4\left( 
 \frac{1}{\sigma (X)} + \frac{1}{\sqrt{\sigma (X)}} \right). 
$$ 
 Margulis-Russo type identities are also obtained and yield asymptotic 
 estimates for the probability of monotone sets. 

 Isoperimetry for graphs and Markov chains is often applied 
 to determine bounds on the spectral gaps $\lambda_2, \lambda_\infty$, 
 providing an estimate of the speed of convergence to equilibrium 
 for stochastic algorithms ans in statistical mechanics. 
 In such situations the values of the isoperimetric constants 
 are easily computed as infima on finite sets. 
 In the configuration space case the situation is different since 
 $\lambda_2$ and $\lambda_\infty$ are known and used to deduce 
 bounds on the isoperimetric constants. 

 We proceed as follows. 
 In Sect.~\ref{s1} we construct a finite difference 
 gradient on Poisson space and recall the associated integration by parts 
 formulas, as well as the Clark formula. We also extend the isoperimetric 
 result of \cite{bobkovgoetze} 
 (see \cite{bobkov2} on 
 Gaussian space and \cite{capitaine} on Wiener space 
 and path space), and state a Margulis-Russo type identity, in the 
 general setting of configuration spaces under Poisson measures. 
 In Sect.~\ref{s2}, 
 a graph is constructed on configuration space by addition 
 or deletion of configuration points. 
 The inner and outer boundaries of subsets of configurations 
 and their surface measures are defined in Sect.~\ref{s3}, 
 e.g. a configuration $\omega \in A$ belongs to the inner boundary 
 of $A$ if it has ``at least'' a (forward or backward) neighbor in $A^c$. 
 A deviation result in terms of the intensity parameter is obtained 
 from the Margulis-Russo identity on Poisson space. 
 Co-area formulas for the finite difference gradient, 
 which differ from the Gauss type formulas of \cite{finkel}, 
 are proved in Sect.~\ref{s4}. 
 Boundary measures and surface measures are defined by averaging 
 the norms of finite difference gradients, 
 which represent the measure of the flow in and out a given set. 
 An equivalence criterion for functional inequalities is also proved. 
 In Sect.~\ref{s7}, the main isoperimetric constants are introduced, 
 and bounds are stated on these constants. 
 Sect.~\ref{s9} is devoted to a generalization of Cheeger's 
 inequality, 
 following the arguments of \cite{bht}, \cite{houdremixed}, \cite{ht}. 
\section{Preliminaries} 
\label{s1} 
 Let $X$ be a metric space with Borel $\sigma$-algebra ${\cal B}(X)$ 
 and let $\sigma$ be a finite and diffuse measure on $X$. 
 Let $\Omega$ denote the set of Radon measures 
$$\Omega = \left\{ 
 \omega = \sum_{i=1}^{i=N} \delta_{x_i} \ : \ 
 (x_i)_{i=1}^{i=N} \subset X, \ x_i\not=x_j, \ \forall i\not= j, \ 
 N\in \inte\cup \{ \infty \}\right\},$$ 
 where $\delta_x$ denotes the Dirac measure at $x\in X$. 
 For convenience of notation 
 we identify $\omega = \sum_{i=1}^{i=n} \delta_{x_i}$ with the 
 set $\omega = \{ x_1,\ldots , x_n \}$. 
 Let ${\cal F}$ denote the $\sigma$-algebra 
 generated by all applications of the form $\omega \mapsto \omega (B)$, 
 $B\in {\cal B}(X)$, and let 
 $\pi$ denote the Poisson measure with intensity $\sigma$ 
 on $\Omega$, 
 defined via 
$$\pi (\{ 
 \omega \in \Omega \ : \  
 \omega (A_1) = k_1, \ldots , \omega (A_n) = k_n 
 \} 
 ) 
 = \prod_{i=1}^{i=n} 
 \frac{\sigma (A_i)^{k_i}}{k_i!} 
 e^{-\sigma (A_i)}
, 
 \quad \quad k_1,\ldots , k_n \in\inte, 
$$ 
 on the $\sigma$-algebra ${\cal F}$ generated by sets of the form 
$$\{ 
 \omega \in \Omega \ : \  
 \omega (A_1) = k_1, \ldots , \omega (A_n) = k_n 
 \},  
$$ 
 for $k_1,\ldots, k_n\in\inte$, 
 and disjoint $A_1,\ldots ,A_n\in {\cal B}(X)$. 
 Let $I_n (f_n)$ denote the multiple Poisson stochastic integral 
 of the symmetric function $f_n\in L^2(X,\sigma )^{\circ n}$, 
 defined as 
$$I_n (f_n) (\omega ) = 
 \int_{\Delta_n} 
 f_n(t_1,\ldots ,t_n) (\omega (dt_1)-\sigma (dt_1)) 
 \cdots (\omega (dt_n)-\sigma (dt_n)), 
 \ \ \ f_n\in L^2_\sigma (X^n)^{\circ n},$$ 
 with 
$\Delta_n = {\{(t_1,\ldots ,t_n)\in X^n \ : \ t_i\not= t_j , \ \forall i\not= j\}}$. 
 We recall the isometry formula 
$$E [I_n(f_n)I_m(g_m)] = n! 
 {\bf 1}_{\{ n=m \}} 
 \langle f_n , g_m\rangle_{L^2_\sigma (X)^{\circ n}},$$ 
 see \cite{nualartvives}. 
 As is well-known, every square-integrable 
 random variable $F\in L^2(\Omega^X, P )$ admits 
 the Wiener-Poisson decomposition 
$$F=\sum_{n=0}^\infty I_n (f_n)$$ 
 in series of multiple stochastic integrals. 
\\ 
 The gradient chosen here on Poisson space is a finite difference operator 
 (see \cite{bpr} for a different construction using derivation operators). 
\begin{definition} 
 For any $F : \Omega \longrightarrow \real$, let 
$$D_x F(\omega ) = 
 (F(\omega ) - F(\omega + \delta_x))1_{\{ 
 x\in \omega^c \}} 
 + (F(\omega ) - F(\omega -\delta_x ))1_{\{x\in \omega \}}, 
$$ 
 for all $\omega\in \Omega$ and $x\in X$. 
\end{definition} 
 Now, given $u : \Omega \times X \to \real$ with sufficient integrability 
 properties, we let 
$$\delta_\sigma (u) 
 = \int_X u(x,\omega ) \sigma (dx) 
 - \int_X u(x,\omega-\delta_x ) \omega (dx), 
$$ 
 and 
$$\delta_\omega (u) 
 = \int_X u(x,\omega ) \omega (dx) 
 - \int_X u(x,\omega + \delta_x ) \sigma (dx) 
. 
$$ 
 Note that in the definition of $\delta_\omega (u)$, the 
 integral over the diffuse measure $\sigma$ makes sense since 
 $\sigma (dx)$-a.s., $x\notin \omega$. 
 Note that 
$$D_x F(\omega + \delta_x )= F(\omega + \delta_x)-F(\omega) 
 = - D_x F(\omega ) , 
 \quad x\notin \omega, 
$$ 
 and 
$$D_x F(\omega - \delta_x )= F(\omega - \delta_x)-F(\omega) 
 = - D_x F (\omega ) , 
 \quad x\in \omega. 
$$ 
 The following relations are then easily obtained: 
\begin{align} 
\label{lll1} 
& \delta_\sigma (uF) 
 = F\delta_\sigma (u) 
 + \delta_\sigma (uDF) 
 - \langle u,DF\rangle_{L^2 (X,\sigma)}, 
\\ 
\label{l2} 
& \delta_\omega (uF) 
 = F\delta_\omega (u) 
 + \delta_\omega (uDF) 
 - \langle u,DF\rangle_{L^2 (X,\omega)}, 
\end{align} 
 and 
\begin{align}
\label{n1} 
& \delta_\sigma (u) 
 = \int_X u(x,\omega ) (\sigma (dx) - \omega (dx)) 
 + \int_X D_x u(x,\omega ) \omega (dx), 
\\ 
\label{n2} 
& \delta_\omega (u) 
 = \int_X u(x,\omega ) (\omega (dx) - \sigma (dx)) 
 + \int_X D_xu(x,\omega ) \sigma (dx). 
\end{align} 
 As shown in Prop.~\ref{adj} below, 
 the operators $\delta_\sigma$ and $\delta_\omega$ 
 are adjoint of $D$, with respect to scalar products 
 respectively given by $\sigma$ and $\omega$. 
\begin{prop} 
\label{adj} 
 We have for $F:\Omega \to \real$ and $v : \Omega \times X \to \real$: 
\begin{equation} 
\label{m1} 
E[F\delta_\sigma (v )] 
 = E[\langle DF, v \rangle_{L^2(\sigma )} 
 ], 
\end{equation} 
 and 
\begin{equation} 
\label{m2} 
E[F\delta_\omega (v )] 
 = E[\langle DF, v \rangle_{L^2(\omega )} 
 ], 
\end{equation} 
 provided the corresponding quantities are integrable. 
\end{prop} 
\begin{Proof} 
 We first show that $E[\delta_\omega (v )] = 0$. 
 For simple processes, this can be proved using 
 the characteristic function of 
 $\int_X h d\omega$ which satisfies 
$$E \left[ \exp \left( i z \int_X h d\omega \right) 
 \right] 
 = 
 \exp \int_X  (e^{izh}-1)d\sigma 
, \quad 
z\in \real 
. 
$$ 
 Differentiating each of those two expressions 
 with respect to $z$ yields 
$$ 
 E \left[ \int_X h d \omega 
 \exp \left( i z \int_X h d\omega \right) 
 \right] 
 = 
 E \left[ 
 \int_X he^{izh} d\sigma 
 \exp \left( i z \int_X h d\omega \right) 
 \right], 
$$ 
 hence 
\begin{eqnarray*} 
 E \left[ \int_X h d(\sigma - \omega ) \exp \left( i z \int_X h d\omega 
 \right) 
 \right] 
& = & 
 E \left[ \langle h,1-e^{izh} \rangle_{L^2(X,\sigma)} 
 \exp \left( i z\int_X h d\omega \right) 
 \right] 
\\ 
& = & 
 E \left[ \left< h,D 
 \exp \left( i z\int_X h d\omega \right) 
 \right>_{L^2(X,\sigma)} 
  \right] 
, 
\end{eqnarray*} 
 where we used the relation $D_x\exp (iz \int_X hd\omega ) = 
 (1-e^{iz h(x)})\exp (iz \int_X hd\omega )$, $\sigma (dx)$-a.e. 
 From \eqref{n2} this implies $E[\delta_\omega (u)] = 0$ 
 for all $u$ of the form 
$$u=\sum_{i=1}^n 1_{A_i} e^{iz_1\omega (B_1)+\cdots + iz_n\omega (B_n)}.$$ 
 By martingale convergence arguments, e.g. as in the proof of Th.~3.4 of 
 \cite{wuls1}, the formula is extended to general $u$. 
 This in turn implies $E[\delta_\sigma (v )]=0$ from 
 \eqref{n1}, and \eqref{m1} using \eqref{lll1}. 
\end{Proof} 
 Note that the relation $E[\delta_\omega (v )] = 0$ 
 can be seen as a consequence of 
 Th.~1 or Cor.~1 in \cite{picard}, 
 and \eqref{m2} follows from \eqref{l2}. 
 We have 
$$ 
 \delta_\sigma D F (\omega ) 
 = \int_X (F(\omega )-F(\omega +\delta_x))\sigma (dx) 
 - \int_X (F(\omega -\delta_x )-F(\omega ))\omega (dx) 
, 
$$ 
 and 
$$ 
 \delta_\omega D F (\omega ) 
 = \int_X (F(\omega )-F(\omega -\delta_x))\omega (dx) 
 - \int_X (F(\omega +\delta_x )-F(\omega ))\sigma (dx)
, 
$$ 
 so that 
\begin{eqnarray} 
\label{soth} 
 \delta_\sigma D F (\omega ) 
 & = & 
 \delta_\omega D F (\omega ) 
 = \int_X D_x F (\omega ) \omega (dx) 
 + \int_X D_x F (\omega ) \sigma (dx)
 \\ 
\nonumber 
& = & 
 (\sigma (X) + \omega (X)) F(\omega ) 
 - \int_X F(\omega +\delta_x )\sigma (dx) 
 - \int_X F(\omega -\delta_x )\omega (dx). 
\end{eqnarray} 
 From the definition of $I_n(f_n)$ it can also be easily shown that 
$$\delta_\sigma DI_n (f_n) = \delta_\omega DI_n (f_n) = n I_n (f_n),$$ 
 cf. e.g. \cite{prirose}. 
 It follows that the spectral gap of $\delta_\sigma D$ is 
 $\lambda_2 = 1$, a fact which is recovered below by a different method. 
 In the sequel we shall uniquely use the operator $\delta_\sigma$, 
 and denote it by $\delta$. 
 Let 
\begin{eqnarray*} 
D^+_x F (\omega ) & = & \max (0 , D_x F (\omega ) ) 
\\ 
& =&  
 (F(\omega ) - F(\omega + \delta_x))^+1_{\{ 
 x\in \omega^c \}} 
 + 
 ( F(\omega ) - F(\omega -\delta_x ))^+1_{\{ x\in \omega \}}. 
\end{eqnarray*} 
 and 
\begin{eqnarray*} 
D^-_x F (\omega ) & = & - \min (0, D_x F (\omega )) 
\\ 
& = &  
 (F(\omega ) - F(\omega + \delta_x))^-1_{\{ 
 x\in \omega^c \}} 
 + 
 (F(\omega ) - F(\omega -\delta_x ))^-1_{\{x\in \omega \}}. 
\end{eqnarray*} 
 We have $D^+_x F = D^-_x (-F)$, 
$$D^+_x F(\omega + \delta_x )= D^-_x F(\omega ) , 
 \qquad 
 D^-_x F(\omega + \delta_x )= D^+_x F(\omega ) , 
 \quad x\notin \omega, 
$$ 
 and 
$$D^+_x F(\omega - \delta_x ) = D^-_x F (\omega ) , 
 \qquad 
 D^-_x F(\omega - \delta_x ) = D^+_x F (\omega ) , 
 \quad x\in \omega, 
$$ 
 which implies 
\begin{equation} 
\label{11} 
 \delta_\sigma (D^+ F)^p (\omega ) 
 = 
 - \delta_\omega (D^- F)^p (\omega ) 
 = \int_X (D^+_x F (\omega ))^p \sigma (dx) 
 - \int_X (D^-_x F (\omega ))^p \omega (dx), 
\end{equation} 
 and 
\begin{equation} 
\label{22} 
 \delta_\sigma (D^- F)^p (\omega ) 
 = 
 - \delta_\omega (D^+ F)^p (\omega ) 
 = \int_X (D^-_x F (\omega ))^p \sigma (dx) 
 - \int_X (D^+_x F (\omega ))^p \omega (dx). 
\end{equation} 
 We also have 
 $\vert D_xF\vert^p = \vert D^+_xF\vert^p + \vert D^-_xF\vert^p$, 
 and 
$$\vert D F (\omega )\vert_{L^p}^p 
 = \vert D^+ F (\omega )\vert_{L^p}^p 
 + \vert D^- F (\omega )\vert_{L^p}^p 
. 
$$ 
\begin{lemma} 
\label{l1} 
 We have 
$$ 
 E[\vert D^+ F \vert_{L^p(\sigma )}^p ]  
 = E[\vert D^- F \vert_{L^p(\omega )}^p ], 
$$ 
 and 
$$ 
 E[\vert D^- F \vert_{L^p(\sigma )}^p]  
 = E[\vert D^+ F \vert_{L^p(\omega )}^p]. 
$$ 
\end{lemma} 
\begin{Proof} 
 Using \eqref{11} and \eqref{22} we have 
$$ 
  E[\vert D^\pm F \vert_{L^p(\sigma )}^p] 
 - E[\vert D^\mp F \vert_{L^p(\omega )}^p]
 = 
 E[\delta_\sigma ((D^\pm F)^p)] 
 = 0 
. 
$$ 
\end{Proof} 
 Similarly, \eqref{soth} will imply 
\begin{equation} 
\label{wlmpl} 
 E\left[\int_X D_x F \sigma (dx) \right]  
 = 
 - E\left[\int_X D_x F \omega (dx) \right]  
.  
\end{equation} 
 In the particular case $F= 1_{\{\omega (A) = k\}}$, 
 Lemma~\ref{l1} simply states the following 
 easily verified equality: 
$$ 
 E[\vert D^+ 1_{\{\omega (A) = k\}} 
 \vert_{L^1(\sigma )}] 
 = \sigma (A) E[1_{\{\omega (A)=k\}} ] 
 = (k+1) E[ 1_{\{\omega (A)=k+1\}} ] 
 = E[\vert D^- 1_{\{\omega (A) = k\}} 
 \vert_{L^1(\omega )}]. 
$$ 
 We also have 
$$E\left[\vert D^+ F \vert_{L^p(\frac{\sigma+\omega}{2} )}^p \right]  
 = E\left[\vert D^- F \vert_{L^p(\frac{\sigma+\omega}{2} )}^p \right] 
 = 
 \frac{1}{2} 
 E\left[\vert D F \vert_{L^p(\sigma )}^p \right]  
 = 
 \frac{1}{2} 
 E\left[\vert D F \vert_{L^p(\omega )}^p\right] 
, 
$$ 
 in particular the Dirichlet forms 
 ${\cal E}_\sigma (F,G)$ and 
 ${\cal E}_\omega (F,G)$ 
 defined as 
$${\cal E}_\sigma (F,F) = \frac{1}{2} 
 E[\vert D F\vert_{L^2 (\sigma )}^2], 
 \quad 
 {\cal E}_\omega (F,F) = \frac{1}{2} 
 E[\vert D F\vert_{L^2 (\omega )}^2] 
$$ 
 coincide: 
$${\cal E}_\sigma (F,F) 
 = 
 {\cal E}_\omega (F,F). 
$$ 
 This result can also be seen as a consequence of the relation 
 $\delta_\sigma D = \delta_\omega D$, or of Prop.~\ref{adj}. 

 The Clark formula given next yields the predictable representation
 of a random variable using the operator $D$.
 Take $X=[0,1]$ and 
 $\sigma$ the Lebesgue measure and let  
$$N_t (\omega ) = N_{[0,t]} (\omega ) 
 = \omega ([0,t]), \quad t\in\real_+, \quad \omega \in \Omega, 
$$ 
 i.e. $(N_t)_{t\in [0,1]}$ is a standard Poisson process 
 under $\pi$. 
\begin{prop}(\cite{privault}, Th.~1) 
\label{cl}
 We have the following Clark formula, 
 for $F\in L^2(\Omega ,\pi)$:
\begin{equation} 
\label{c1} 
F = E[F] - \int_0^1 E[D_t F \mid {\cal F}_t]d\tilde{N}_t, 
\end{equation} 
 where the stochastic integral is taken in the It\^o sense. 
\end{prop} 
 The formula is first proved for $F\in \Dom (D)$ and then 
 extended to $L^2(\Omega ) $ 
 by continuity of $F\mapsto (E[D_t F \mid {\cal F}_t ])_{t\in \real_+}$ 
 from $L^2(\Omega ,\pi)$ into $L^2(\Omega \times [0,1])$. 
 The Clark formula \eqref{cl} yields the Poincar\'e inequality: 
\begin{equation} 
\label{ll1} 
\var (F) \leq E[\vert DF\vert_{L^2(\sigma )}^2], \quad 
 F\in \Dom (D). 
\end{equation} 
 This inequality is in fact valid for an arbitrary Polish space 
 $X$ with diffuse measure $\sigma$. 
 Note that if $F=1_A$ then the Poincar\'e inequality implies 
$$\pi (A) (1-\pi (A)) \leq 
 \sigma (X),$$ 
 in particular if $\sigma (X) \leq 1/4$ then 
 we have either 
$$\pi (A) \leq (1-\sqrt{1-4\sigma (X)})/2$$ 
 or 
$$\pi (A) \geq (1+\sqrt{1-4\sigma (X)})/2,$$ 
 and if $\pi (A)\leq 1/2$ then 
$$\pi (A) \leq 2 \pi (A) (1-\pi (A)) \leq 
 2 \sigma (X).$$ 
 The following result gives a version of isoperimetry on Poisson space 
 which is independent of dimension and generalizes the result of \cite{bobkovgoetze}, p. 274. 
 Let $\varphi$ denote the standard Gaussian density, and 
 let $\Phi$ denote its distribution function. 
 Let $I(t) = \varphi (\Phi^{-1}(t))$, 
 $0\leq t \leq 1$ denote the Gaussian isoperimetric function, 
 with the relations $I(x)I''(x) = -1$
 and $I'(x) = -\Phi^{-1}(x)$, $x\in [0,1]$. 
\begin{prop}
\label{9.4} 
 For every random variable $F:\Omega \rightarrow [0,1]$ we have 
\begin{equation} 
\label{e1} 
I(E[F]) \leq E\left[\sqrt{I( F )^2 + 2 \vert DF\vert_{L^2 (\sigma )}^2}\right]. 
\end{equation}
\end{prop} 
\begin{Proof} 
 Let $X_n$ denote the $\inte^n$-valued random variable defined as 
$$X_n (\omega ) 
 = \left(\omega (A_1), \ldots , \omega (A_n ) \right), \quad \omega 
 \in \Omega. 
$$ 
 If $F = f \circ X_n$ is a cylindrical functional we have 
$$ 
D_x F (\omega ) = \sum_{k=1}^{k=n} 
 1_{A_k}(x) 
 ( f ( X_n(\omega )) - f(X_n(\omega)+e_k) ), 
$$ 
 $f:\inte^n \to \real$, 
 where $(e_k)_{1\leq k\leq n}$ denotes the canonical basis 
 of $\real^n$. 
 For the cylindrical functional $F$, 
 \eqref{e1} follows by application of Relation (3.13) 
 in \cite{bobkovgoetze} and tensorization. 
 The extension to general random variables can be done by 
 martingale convergence, e.g. as in the proof of Th.~3.4 of 
\cite{wuls1}. 
\end{Proof} 
 This also implies that the optimal constant $b_2$ in the 
 inequality 
$$ 
I(E[F]) \leq E\left[\sqrt{I( F )^2 + \frac{1}{b_2} 
 \vert DF\vert_{L^2 (\frac{\sigma +\omega}{2} )}^2}\right] 
$$ 
 satisfies $b_2\geq 1$. 
 Using the equivalence 
 $I(\varepsilon ) \simeq \varepsilon \sqrt{2\log 1/\varepsilon}$ 
 and the Schwarz inequality, Relation \eqref{e1} 
 allows to recover the modified logarithmic Sobolev 
 inequality of \cite{ane}, \cite{wuls2}: 
$$E[F\log F] - E[F]\log E[F] \leq \frac{1}{2} 
 E\left[ \frac{1}{F} \vert DF\vert_2^2 \right]. 
$$ 
 Note that the analog Gaussian isoperimetry result can also be 
 transferred to the Poisson space for 
 the Carlen-Pardoux gradient \cite{carlen}, writing the exponential 
 interjump times of the Poisson process as half sums of squared 
 Gaussian random variables as in \cite{cras3}. 
 Let $\pi_\lambda$, $\lambda>0$, 
 denote the Poisson measure of intensity $\lambda \sigma 
 (dx)$ on $\Omega$, and let $E_\lambda$ denote the expectation under 
 $\pi_\lambda$. 
 We refer to \cite{zuyev} for the following type of result, obtained 
 by differentiation of the intensity parameter. 
\begin{prop} 
\label{pr2.6} 
 Assume that $DF \in L^1(\pi_\lambda \otimes \sigma )$ 
 and $F \in L^1(\pi_\lambda )$, $\lambda \in (a,b)$. 
 We have 
$$\frac{\partial}{\partial \lambda} 
 E_\lambda [F] 
 = - E_\lambda \left[ \int_X D_x F \sigma (dx)\right] 
 = E_\lambda \left[ \int_X D_x F \omega (dx)\right], 
 \quad \lambda \in (a,b). 
$$ 
\end{prop} 
\begin{Proof} 
 Given the representation 
$$F(\omega ) = f_0 1_{\{\vert \omega \vert = 0 \}} 
 + \sum_{n=1}^\infty 
 1_{\{\vert \omega \vert = n\}} 
 f_n(x_1,\ldots , x_n), 
$$ 
 where $\omega = \{x_1,\ldots ,x_n\}$ when $\vert \omega \vert = n$, 
 we have 
$$E_\lambda [F] 
 = 
 e^{-\lambda \sigma (X)} 
 f_0 
 + 
 e^{-\lambda \sigma (X)} 
 \sum_{n=1}^\infty 
 \frac{\lambda^n}{n!} 
 \int_X \cdots \int_X 
 f_n(x_1,\ldots , x_n) 
 \sigma (dx_1) \cdots \sigma (dx_n), 
$$ 
 and 
\begin{eqnarray*} 
\frac{\partial}{\partial \lambda} 
 E_\lambda [F] 
 & = & 
 -\sigma (X) 
 E_\lambda [F] 
\\ 
& & + 
 e^{-\lambda \sigma (X)} 
 \sum_{n=1}^\infty 
 \frac{\lambda^{n-1}}{(n-1)!} 
 \int_X \cdots \int_X 
 f_n(x_1,\ldots , x_n) 
 \sigma (dx_1) \cdots \sigma (dx_n)
\\ 
& = & 
 -\sigma (X) 
 E_\lambda [F] 
 + 
 E_\lambda \left[ 
 \int_X  F(\omega + \delta_x) \sigma (dx) 
 \right] 
\\ 
& = & 
 - E_\lambda \left[ 
 \int_X  D_xF(\omega ) \sigma (dx) 
 \right] 
. 
\end{eqnarray*} 
 The second relation follows from \eqref{wlmpl}. 
\end{Proof} 
 As a corollary we will obtain a Margulis-Russo type equality 
 \cite{margulis}, \cite{lrusso} for monotone sets under Poisson measures. 
\begin{definition} 
 A measurable set $A\subset \Omega$ is called increasing 
 if 
\begin{equation} 
\label{propr1} 
\omega \in A \quad \Longrightarrow \quad 
 \omega + \delta_x \in A, \quad \sigma (dx)-a.e. 
\end{equation} 
 It is called decreasing if 
\begin{equation} 
\label{propr2} 
\omega \in A \quad \Longrightarrow \quad 
 \omega - \delta_x \in A, \quad \omega (dx)-a.e. 
\end{equation} 
\end{definition} 
 Note that if $A$ is decreasing then $A^c$ is increasing but the converse 
 is not true. 
 In fact, saying that $A$ is decreasing is equivalent to 
 the following property on $A^c$: 
\begin{equation} 
\label{propr} 
\omega \in A^c \quad \Longrightarrow \quad \omega + \delta_x \in A^c, \quad 
 \forall x\in \omega^c, 
\end{equation} 
 which is stronger than saying that $A^c$ is increasing. 
 The set $A$ is said to be monotone if it is either increasing or decreasing. 
 The sets $\{ \omega (B) \geq n\}$, resp. 
 $\{ \omega (B) \leq n\}$, are naturally increasing, resp. decreasing. 
 Another example of monotone set is given by 
$$\left\{ 
 \omega \in \Omega \ : \  
 \int_X fd\omega >K \right\}, \quad K\in \real, 
$$ 
 which is increasing, resp. decreasing, if $f\geq 0$, resp. $f\leq 0$. 
 Clearly, a set $A$ is increasing, resp. decreasing, if and only if 
 $D_x1_A\leq 0$ (i.e. $D_x1_A = - D^-_x1_A$, or $D^+_x 1_A =0$) 
 $\sigma (dx)$-a.e., resp. $\omega (dx)$-a.e. 
 As a corollary of Prop.~\ref{pr2.6} we have: 
\begin{corollary} 
\label{margulisrusso} 
 Let $A\subset \Omega$ be an increasing set. 
 We have 
$$\frac{\partial}{\partial \lambda} 
 \pi_\lambda (A) 
 = E_\lambda \left[ \int_X D^-_x 1_A \sigma (dx)\right] 
 = E_\lambda \left[ \int_X D^+_x 1_A \omega (dx)\right] 
. 
$$ 
 If $A\subset \Omega$ is decreasing we have 
$$\frac{\partial}{\partial \lambda} 
 \pi_\lambda (A) 
 = - E_\lambda \left[ \int_X D^-_x 1_A \omega (dx)\right] 
 = - E_\lambda \left[ \int_X D^+_x 1_A \sigma (dx)\right] 
. 
$$ 
\end{corollary} 
 We also have if $A$ is monotone: 
$$\frac{\partial}{\partial \lambda} 
 \pi_\lambda (A) 
 = E_\lambda \left[ \Vert D 1_A\Vert_{L^1(\sigma )} \right] 
 = E_\lambda \left[ \Vert D 1_A\Vert_{L^1(\omega )} \right]. 
$$ 
\section{Forward-backward kernels and reversibility on configuration space} 
\label{s2} 
 Given $\omega \in \Omega$, the set of forward neighbors of 
 $\omega$ is defined to be 
$${\cal N}^{+}_\omega = \{ \omega + \delta_x \ : x\in \omega^c \}, 
$$ 
 and similarly 
 the set of backward neighbors of 
 $\omega$ is 
$${\cal N}^{-}_\omega = \{ \omega - \delta_x \ : x\in \omega \}. 
$$ 
 We let 
$${\cal N}_\omega  = {\cal N}^{+}_\omega \cup {\cal N}^{-}_\omega. 
$$ 
 We define two measure kernels $K^+(\omega , d\tilde{\omega})$ 
 and $K^-(d\tilde{\omega} , \omega )$ which are respectively 
 supported by ${\cal N}_\omega^+$ and ${\cal N}_\omega^-$. 
\begin{definition} 
 Let for $A\in {\cal F}$: 
$$K^+(\omega , A ) 
 = \int_X 1_A (\omega +\delta_x)\sigma (dx), 
\quad \quad 
K^-(A , \omega ) = 
 \sum_{x \in \omega } 1_A (\omega-\delta_x). 
$$ 
\end{definition} 
 It is a classical fact that since $\pi$ is a Poisson measure, 
 the image under $\omega+\delta_x\mapsto x$ 
 of the measure 
$$\pi (d\tilde{\omega} \mid \tilde{\omega} \in {\cal N}^{+}_\omega 
 ) 
$$ 
 coincides with the (normalized) measure $\sigma$ on $X$: 
$$
 \frac{\sigma (B)}{\sigma (X)} 
 = 
 \pi ( 
 \{ \tilde{\omega} \ : \ \tilde{\omega} = \omega + \delta_x \ : x\in B\} 
 \mid \tilde{\omega} \in {\cal N}^{+}_\omega  
 ), \quad 
 B\in {\cal B}(X). 
$$ 
 Hence the forward kernel satisfies 
$$K^+(\omega, d\tilde{\omega} ) = \sigma (X) 
 \pi (d\tilde{\omega} \mid \tilde{\omega} \in {\cal N}^{+}_\omega 
 ), 
$$ 
 and $(\sigma (X))^{-1} K^+(\omega , d\tilde{\omega } )$ is of Markov type. 
 Similarly, the image under $\omega-\delta_x \mapsto x$ 
 of the measure 
$$\pi (d\tilde{\omega} \mid \tilde{\omega} \in {\cal N}^{-}_\omega 
 ) 
$$ 
 coincides with the normalized counting measure on $\omega$: 
$$
\frac{\omega (B)}{\omega (X)} 
 =  
 \pi (
 \{ \tilde{\omega} \ : \  
 \tilde{\omega} = \omega - \delta_x, \ x \in B\} 
 \mid \tilde{\omega} \in {\cal N}^{-}_\omega  
 ), 
$$ 
 hence the backward kernel satisfies 
$$K^-(d\tilde{\omega} , \omega ) = 
 \omega (X) 
 \pi (d\tilde{\omega} \mid \tilde{\omega} \in {\cal N}^{-}_\omega 
 ) 
 = \sum_{x \in \omega } \delta_{\omega-\delta_x} (d\tilde{\omega}), 
$$ 
 and 
$(\omega (X))^{-1}K^-(d\tilde{\omega},\omega )$ is Markovian 
 provided $\omega \not= \emptyset$. 
 The kernel $K^-(d\tilde{\omega},\omega )$ 
 itself is not Markovian, instead it is of graph type, i.e. 
$$K^-(\{ \tilde{\omega}\}, \omega ) = 
 \left\{ 
 \begin{array}{ll} 
 1 & \mbox{if } \tilde{\omega} = \omega - \delta_x \mbox{ for some } x\in X 
 \ (\mbox{i.e. } \tilde{\omega}\in {\cal N}_\omega ), 
\\ 
 0 & \mbox{otherwise} 
 \ (\mbox{i.e. } \tilde{\omega}\notin {\cal N}_\omega ). 
\end{array} 
\right. 
$$ 
 We have for $p\in [1,\infty )$: 
$$\vert D F (\omega )\vert_{L^p(\sigma )}^p = 
 \int_X \vert F(\omega ) - F(\omega + \delta_x )\vert^p 
 \sigma (dx) 
 = 
 \int_\Omega \vert F(\omega ) - F(\tilde{\omega} )\vert^p K^+(\omega , 
 d\tilde{\omega}), 
$$ 
 and 
$$\vert D F (\omega )\vert_{L^p(\omega )}^p = 
 \int_X \vert F(\omega ) - F(\omega -\delta_x )\vert^p 
 \omega (dx ) 
 = 
 \int_\Omega \vert F(\omega ) - F(\tilde{\omega} )\vert^p K^-(\omega , 
 d\tilde{\omega}). 
$$ 
 For $p = \infty $ we have 
$$
\vert D F (\omega )\vert_{L^\infty (\sigma )} = 
 \esssup_{\sigma (dx)} \vert F(\omega ) - F(\omega + \delta_x )\vert 
 = 
 \esssup_{K^+(\omega, d\tilde{\omega})} 
 \vert F(\omega ) - F(\tilde{\omega} )\vert 
, 
$$ 
 and 
$$ 
\vert D F (\omega )\vert_{L^\infty (\omega )} = 
 \esssup_{\omega (dx) } \vert F(\omega ) - F(\omega - \delta_x )\vert 
 = 
 \esssup_{K^-(\omega, d\tilde{\omega})} 
 \vert F(\omega ) - F(\tilde{\omega} )\vert 
. 
$$ 
 We also have 
$$ 
E\left[ \vert D^+ 1_A\vert_{L^p (\frac{\sigma+\omega}{2} )} \right] 
 = \int_A \bar{K}(\omega , A^c)^{1/p} 
 \pi (d\omega ), 
 \quad 
 E\left[ \vert D^- 1_A\vert_{L^p (\frac{\sigma+\omega}{2} )} \right] 
 = \int_{A^c} \bar{K}(\omega , A)^{1/p} 
 \pi (d\omega ). 
$$ 
 The following proposition shows a reversibility property, 
 which is an analog of Lemma~\ref{l1}. 
\begin{prop} 
\label{rv} 
 The kernels $K^+ (\omega , d\tilde{\omega })$ 
 and $K^-(d\tilde{\omega } , \omega )$ are 
 mutually adjoint under $\pi (d\tilde{\omega} )$, i.e. 
$$\pi (d\omega ) K^+(\omega , d\tilde{\omega }) 
 = K^-(d\omega , \tilde{\omega} ) \pi (d\tilde{\omega} ). 
$$ 
\end{prop} 
\begin{Proof} 
 We have 
\begin{eqnarray*} 
 \int_\Omega \int_\Omega F(\omega ) G(\tilde{\omega}) 
 K^+(\omega , d\tilde{\omega }) \pi (d\omega ) 
 & = & 
 \int_\Omega F(\omega ) G(\omega + \delta_x ) 
 \pi (d\omega ) \sigma (dx) 
\\ 
 & = & 
 - E[F \langle DG ,1\rangle_{L^2(\sigma )} ] 
 + \sigma ( X ) E[FG] 
\\ 
 & = & 
 - E[G\delta_\sigma (1_X F) ] 
 + \sigma ( X ) E[FG] 
\\ 
 & = & 
 \int_\Omega G(\omega ) \sum_{x\in \omega} 
 F(\omega - \delta_x ) 
 \pi (d\omega ) 
\\ 
 & = & 
 \int_\Omega \int_\Omega G(\tilde{\omega} ) F( \omega ) 
 K^- (d\omega , \tilde{\omega} ) \pi (d\tilde{\omega} ) . 
\end{eqnarray*} 
\end{Proof} 
 In particular we have $E[K^- F] = \sigma (X) E[F]$: 
$$\int_\Omega \int_\Omega 
 F(\tilde{\omega}) K^-(d\tilde{\omega},\omega ) 
 \pi (d\omega ) 
 = \sigma (X) \int_\Omega F(\omega ) \pi (d\omega ), 
$$ 
 and $E[K^+ F] =E[\omega (X)F]$: 
$$\int_\Omega \int_\Omega 
 F(\tilde{\omega}) K^+(\omega , d\tilde{\omega}) 
 \pi (d\omega ) 
 = \int_\Omega \omega (X) F(\omega ) \pi (d\omega ), 
$$ 
 which is Lemma~1.1 in \cite{wuls1} and is similar to the Mecke identity 
 \cite{jmecke}. 
 This also implies 
\begin{align*} 
& \int_{A} K^+(\omega ,A^c) \pi (d\omega ) 
 = E[\vert D^+ 1_A (\omega ) \vert_{L^p(\sigma )}^p] 
 = E[\vert D^- 1_A (\omega ) \vert_{L^p(\omega )}^p] 
 = \int_{A^c} K^-(A,\omega ) \pi (d\omega ), 
\\ 
& \int_{A^c} K^+(\omega ,A) \pi (d\omega ) 
 = E[\vert D^- 1_A (\omega ) \vert_{L^p(\sigma )}^p] 
 = E[\vert D^+ 1_A (\omega ) \vert_{L^p(\omega )}^p] 
 = \int_{A} K^-(A^c,\omega ) \pi (d\omega ) 
. 
\end{align*} 
 The proof of Lemma~\ref{l1} can be reformulated using reversibility 
 of forward and backward kernels. 
\begin{Proof} 
 We have 
\begin{eqnarray*} 
E[\vert D^\pm F \vert_{L^p(\sigma )}^p ]  
 & = & 
\int_{\Omega} 
 ((F(\omega ) - F(\tilde{\omega} ))^\pm)^p 
 K^+(\omega , d\tilde{\omega} ) \pi (d\omega ) 
\\ 
 & = & 
\int_{\Omega} 
 ((F(\omega ) - F(\tilde{\omega} ))^\pm)^p 
 K^-(d\omega , \tilde{\omega} ) \pi (d\tilde{\omega} ) 
\\ 
 & = & 
\int_{\Omega} 
 ((F(\tilde{\omega} ) - F(\omega ))^\mp)^p 
 K^-(d\omega , \tilde{\omega} ) \pi (d\tilde{\omega} ) 
\\ 
& = & E[\vert D^\mp F \vert_{L^p(\omega )}^p ]. 
\end{eqnarray*} 
\end{Proof} 
\noindent 
 Let $\bar{K}(\omega , d\tilde{\omega})$ 
 denote the symmetrized kernel 
$$\bar{K}(\omega , d\tilde{\omega}) 
 = \frac{K^+(\omega , d\tilde{\omega}) 
 + K^-(d\tilde{\omega} , \omega )}{2}. 
$$ 
 We have 
$$\vert D F (\omega )\vert_{L^p(\frac{\omega +\sigma}{2})}^p = 
 \frac{1}{2} 
 \vert D F (\omega )\vert_{L^p(\sigma )}^p 
 + 
 \frac{1}{2} 
 \vert D F (\omega )\vert_{L^p(\omega )}^p 
 = 
 \int_\Omega \vert F(\omega ) - F(\tilde{\omega} )\vert^p \bar{K} (\omega , 
 d\tilde{\omega}), 
$$ 
 and for $p = \infty $: 
$$
\vert D F (\omega )\vert_{L^\infty (\sigma + \omega)} = 
 \esssup_{\bar{K} (\omega, d\tilde{\omega})} 
 \vert F(\omega ) - F(\tilde{\omega} )\vert 
. 
$$ 
 We also have 
\begin{eqnarray*} 
E\left[ \vert D1_A\vert_{L^p (\frac{\sigma+\omega}{2} )} \right] 
& = & E\left[ \vert D^+ 1_A\vert_{L^p (\frac{\sigma+\omega}{2} )} \right] 
 + E\left[ \vert D^- 1_A\vert_{L^p (\frac{\sigma+\omega}{2} )} \right] 
\\ 
&  = & \int_A \bar{K}(\omega , A^c)^{1/p} 
 \pi (d\omega ) 
 + 
 \int_{A^c} \bar{K}(\omega , A)^{1/p} 
 \pi (d\omega ), 
\end{eqnarray*} 
 since $D^+_x F D^-_x F =0$, $x\in X$. 
 Let 
$$\Gamma^\pm (F,F) 
 = \frac{1}{2} 
 \vert D^\pm F\vert_{L^2 (\frac{\sigma+\omega}{2} )}^2. 
$$ 
 We have 
$$\Gamma^+ (F,G)(\omega )  
 = \frac{1}{2} 
 \int_{\Omega\times \Omega} 
 (F(\omega ) - F(\tilde{\omega} ) )^+ 
 (G(\omega ) - G(\tilde{\omega} ) )^+ 
 \bar{K}(\omega , d\tilde{\omega} ),  
$$ 
$$\Gamma^- (F,G)(\omega )  
 = \frac{1}{2} 
 \int_{\Omega\times \Omega} 
 (F(\omega ) - F(\tilde{\omega} ) )^- 
 (G(\omega ) - G(\tilde{\omega} ) )^- 
  \bar{K}(\omega , d\tilde{\omega} ),  
$$ 
 and 
$$ 
{\cal E}(F,F) 
= E[\Gamma^+ (F,F)] 
 = E[\Gamma^- (F,F)]
. 
$$ 
\begin{prop} 
 The Laplacian associated to the discrete Dirichlet form 
 ${\cal E}(F,F)$ is $L = \frac{1}{2}\delta D$, 
 with 
$$L = \frac{1}{2} 
 \delta D =  \frac{\sigma (X) + \omega (X)}{2} I_d - \bar{K}. 
$$ 
\end{prop} 
\begin{Proof} 
 Again, reversibility can be employed. 
 We have 
\begin{eqnarray*} 
 {\cal E}(F,G) 
 & = & \int_{\Omega \times \Omega} 
 (F(\omega )- F(\tilde{\omega} )) 
 (G(\omega )- G(\tilde{\omega} )) 
 K^+(\omega , d\tilde{\omega} ) 
 \pi (d\omega ) 
\\ 
 & = & 
 \int_{\Omega \times \Omega} 
 F(\omega ) G(\omega ) 
 K^+(\omega , d \tilde{\omega})  
 \pi (d\omega ) 
 + 
 \int_{\Omega \times \Omega} 
 F(\tilde{\omega} ) G(\tilde{\omega} ) 
 K^-(d\omega , \tilde{\omega})  
 \pi (d\tilde{\omega} ) 
\\ 
& & - \int_{\Omega \times \Omega} 
 F(\omega ) G(\tilde{\omega} ) 
 K^+(\omega , d\tilde{\omega} ) 
 \pi (d\omega ) 
 - 
 \int_{\Omega \times \Omega} 
 G(\omega ) F(\tilde{\omega} ) 
 K^+(\omega , d\tilde{\omega} ) 
 \pi (d\omega ) 
\\ 
 & = & 
 E[F((\sigma (X) + \omega (X)) G 
 - K^+G-K^-G) ]. 
\end{eqnarray*}  
\end{Proof} 
 Note that in the case of cylindrical functionals, 
 $L$ is the generator of Glauber dynamics considered in statistical 
 mechanics as in e.g. \cite{daipra}, and has the Poisson 
 probability as invariant measure. 
 Although $K^-(d\tilde{\omega},\omega)$ 
 and $K^+(\omega , d\tilde{\omega})$ 
 are not Markov, they leave the Poisson measure invariant 
 under appropriate normalizations, for example 
 for $A=\{ \omega ( X ) = k\}$, we have 
 $K^-(A,\omega ) = (k+1) 
 1_{\{\omega (X) = k+1\}}$, 
 and 
$$\frac{1}{\sigma (X)} 
 \int_\Omega \pi (d\omega ) K^-(A,\omega ) 
 = \frac{k+1}{\sigma (X)} 
 \pi (\{ \omega (X) = k+1 \}) 
 = \pi ( A ). 
$$ 
 In particular we have the following result. 
\begin{prop} 
 The Poisson measure $\pi (d\omega )$ 
 is a stationary distribution for the symmetrized normalized 
 kernel 
$$
 \frac{2}{\sigma (X) 
 + \tilde{\omega} (X)} \bar{K}(\omega ,d\tilde{\omega})
. 
$$ 
\end{prop} 
\begin{Proof} 
 We have 
$$\int_\Omega \pi (d\omega ) 
 \frac{2}{\sigma (X)+\omega (X)} \bar{K}(\omega , A ) 
 = 
 \int_A \pi (d\omega ) 
 \frac{2}{\sigma (X)+\omega (X)} 
 \bar{K}(\omega , \Omega ) 
 = \pi (A). 
$$ 
\end{Proof} 
\section{Inner and outer boundaries} 
\label{s3} 
 We have 
$$D^+_x 1_A (\omega ) = 
 1_{\{\omega \in A \ \mbox{and} \  
 \omega +\delta_x \in A^c\}}1_{\{ 
 x\in \omega^c \}} 
 + 
 1_{\{\omega \in A \ \mbox{and} \  
 \omega - \delta_x \in A^c \}} 1_{\{x\in \omega \}}, 
$$ 
 and 
$$D^-_x 1_A (\omega ) = 
 1_{\{\omega \in A^c \ \mbox{and} \  
 \omega +\delta_x \in A\}} 
 1_{\{ x\in \omega^c \}} 
 + 
 1_{\{\omega \in A^c \ \mbox{and} \  
 \omega - \delta_x \in A \}} 1_{\{ 
 x\in \omega \}}. 
$$ 
 Hence 
$$\vert D^+ 1_A (\omega ) \vert_{L^p(\sigma )}^p 
 = 1_A (\omega ) \sigma (\{ 
 x \in X \ : \  
 \omega +\delta_x \in A^c\}) 
 = 1_A (\omega ) K^+(\omega , A^c), 
$$ 
 and 
$$\vert D^+ 1_A (\omega ) \vert_{L^p(\omega )}^p 
 = 
 1_{A} (\omega ) \omega (\{ 
 x \in X \ : \  
 \omega - \delta_x \in A^c \}) 
 = 1_A(\omega ) K^-(A^c,w), 
$$ 
 i.e. 
 for $\omega\in A$, 
 $\vert D^+ 1_A (\omega ) \vert_{L^p(\sigma )}^p$ is the measure 
 $K^+(\omega , A^c)$ on ${\cal N}^+_\omega$ 
 of the set of forward neighbors which belong to $A^c$, 
 and 
 $\vert D^+ 1_A (\omega ) \vert_{L^p(\omega )}^p$ is the 
 number (or measure $K^-(A^c ,\omega )$ on ${\cal N}^-_\omega$) 
 of backward neighbors which belong to $A^c$. 
 We also have 
$$\vert D^- 1_A (\omega ) \vert_{L^p(\sigma )}^p  
 = 1_{A^c} (\omega ) 
 \sigma (\{ 
 x \in [0,1] \ : \  
 \omega +\delta_x \in A\}) 
 = 1_{A^c} (\omega ) K^+(\omega , A), 
$$ 
 and 
$$\vert D^- 1_A (\omega ) \vert_{L^p (\omega )}^p  
 = 1_{A^c} (\omega ) 
 \omega (\{ x \in [0,1] \ : \  
 \omega -\delta_x \in A\}) 
 = 1_{A^c}(\omega ) K^-(A,\omega ). 
$$ 
 i.e. for $\omega\in A^c$, 
 $\vert D^- 1_A (\omega ) \vert_{L^p(\sigma )}^p$ is the measure 
 $K^+(\omega , A )$ on ${\cal N}^+_\omega$ 
 of the set of forward neighbors of $\omega\in A^c$ which belong 
 to $A$, and $\vert D^- 1_A (\omega ) \vert_{L^p(\omega )}^p$ 
 is the number (measure $K^-(A,\omega )$ on ${\cal N}^-_\omega$) 
 of backward neighbors of $\omega\in A^c$ which belong 
 to $A$. 
\begin{remark} 
\label{r1} 
 We have $D^+_x 1_A = D^-_x 1_{A^c}$ and $\vert D_x1_A\vert 
 = \vert D_x1_{A^c}\vert$, $x\in X$.  
\end{remark} 
 In particular, 
$$D^+_x 1_{\{ \omega (B) = k\}} = 
 1_B (x) 1_{\{ \omega (B) = k\}} 
, 
$$ 
 and 
$$D^-_x 1_{\{ \omega (B) = k\}} = 
 1_B (x) 1_\omega (x) 1_{\{ \omega (B) = k+1\}} 
 + 
 1_B (x) 1_{\omega^c} (x) 1_{\{ \omega (B) = k-1\}} 
, 
$$ 
 hence 
$$\vert D^+ 1_{\{ \omega (B) = k\}} \vert_{L^p(\sigma )}^p 
 = \sigma (B) 1_{\{ \omega (B) = k\}} 
, 
 \quad 
 \vert D^+ 1_B (\omega ) \vert_{L^p(\omega )}^p 
 = 
 k 1_{\{ \omega (B) = k\}} 
, 
$$ 
 and 
$$\vert D^- 1_{\{ \omega (B) = k\}} \vert_{L^p(\sigma )}^p  
 = 
\sigma (B) 1_{\{ \omega (B) = k-1\}} 
, 
  \quad 
\vert D^- 1_{\{ \omega (B) = k\}} \vert_{L^p (\omega )}^p  
 = (k+1) 1_{\{ \omega (B) = k+1\}} 
. 
$$ 
 Similarly, 
\begin{eqnarray*} 
& & \vert D^+ 1_A (\omega ) \vert_{L^\infty (\sigma )} 
 = 
 1_{\{\omega \in A \ \mbox{and} \ 
 \sigma (\{ 
 x \in X \ : \  
 \omega +\delta_x \in A^c\}) >0\}} 
 = 1_A (\omega ) 1_{\{K^+(\omega , A^c)>0\}} 
 , 
\\ 
& & \vert D^+ 1_A (\omega ) \vert_{L^\infty (\omega ) } 
 = 1_{\{\omega \in A \ \mbox{and} \ 
 \exists x \in \omega \ : \  
 \omega - \delta_x \in A^c \}} 
 = 1_A (\omega ) 1_{\{K^-(A^c, \omega )>0\}}, 
\\ 
& & \vert D^- 1_A (\omega ) \vert_{L^\infty (\sigma )} 
 = 
 1_{\{\omega \in A^c \ \mbox{and} \ 
 \sigma (\{ x\in X \ : \  
 \omega +\delta_x \in A\}) >0 \}} 
 = 1_{A^c} (\omega ) 1_{\{K^+(\omega , A )>0\}} 
 ,  
\\ 
& & \vert D^- 1_A (\omega ) \vert_{L^\infty (\omega )} 
 = 
 1_{\{\omega \in A^c \ \mbox{and} \  
 \exists x \in \omega \ : \  
 \omega - \delta_x \in A\}} 
  = 1_{A^c} (\omega ) 1_{\{K^-(A, \omega )>0\}}, 
\end{eqnarray*} 
 i.e. $\vert D^+ 1_A (\omega ) \vert_{L^\infty (\sigma )} = 1$, 
 resp. $\vert D^- 1_A (\omega ) \vert_{L^\infty (\sigma )} =1$, if and only 
 if $\omega\in A$, resp. $\omega\in A^c$, 
 has ``at least'' a forward neighbor in $A^c$, resp. $A$, 
 and $\vert D^+ 1_A (\omega ) \vert_{L^\infty (\omega )} = 1$, 
 resp. $\vert D^- 1_A (\omega ) \vert_{L^\infty (\omega )} =1$, if and only 
 if $\omega\in A$, resp. $\omega\in A^c$, 
 has at least a backward neighbor in $A^c$, resp. $A$. 
 The following definitions are stated independently of $p\in [1,\infty]$. 
\begin{definition} 
 Let $p\in [1,\infty ]$. 
\begin{description} 
\item The inner and outer boundaries of $A$ are defined as: 
$$\partial_{\rm in} A 
 = \{ \omega \in A \ : \ 
 \bar{K}(\omega,A^c) >0\} 
 = \{ 
 \vert D^+ 1_A (\omega ) \vert_{L^p (\sigma+\omega )} >0 \}  
, 
$$ 
 and 
$$\partial_{\rm out} A = \{ \omega \in A^c \ : \ 
 \bar{K}(\omega,A) >0\} 
 = \{ 
 \vert D^- 1_A (\omega ) \vert_{L^p (\sigma+\omega )} >0 \} 
. 
$$ 
\end{description} 
\end{definition} 
 The boundary of $A$ is defined as: 
\begin{eqnarray*} 
 \partial A & = & \partial_{\rm in} A\cup \partial_{\rm out} A 
\\ 
& = & \{ \omega \in \Omega \ : \  
 \vert D 1_A (\omega ) \vert_{L^p (\sigma+\omega )} >0 \} 
\\ 
& = & \{ \omega \in \Omega \ : \ 
 \bar{K}(\omega,A)+\bar{K}(\omega,A^c) >0\}. 
\end{eqnarray*} 
 For instance, 
\begin{align*} 
&\partial_{\rm in} \{ \omega (B ) = k \} 
 = \{ \omega (B ) = k  \}, 
\\ 
\\ 
& \partial_{\rm out} \{ \omega (B ) = k \} 
 = \{ \omega (B ) = k-1 \} \cup 
 \{ \omega (B ) = k+1 \} , 
\\ 
\\
&\partial \{ \omega (B ) = k \} 
 = \{ k-1\leq \omega (B ) \leq k+1  \}. 
\end{align*} 
 In particular, Prop.~\ref{9.4} shows that 
 the isoperimetric function $p\mapsto 
 \inf_{\pi (A) = p} \pi_s (\partial A)$ 
 on Poisson space is greater than $1/\sqrt{2}$ times the 
 Gaussian isoperimetric function $I$. 
 We have $D^+_x 1_A = D^-_x 1_{A^c}$, hence 
 $\partial_{\rm in} A = \partial_{\rm out} A^c$ 
 and $\partial A = \partial A^c$. 
 We may also define the interior $A^\circ$ of $A$ as 
$$A^\circ = \{ \omega \ : \ 
\vert D^+ 1_A (\omega ) \vert_{L^p(\frac{\sigma+\omega}{2} )} = 0 \} 
 = 
\{ \omega \in A \ : \ 
 \bar{K} (\omega , A ) = 0 \} 
 = A\setminus \partial_{\rm in} A,$$ 
 and the closure $\bar{A}$ of $A$ as 
\begin{eqnarray*} 
 \bar{A} & = & \{ \omega \in A^c \ : \ 
 \vert D^- 1_A (\omega ) \vert_{L^p(\frac{\sigma+\omega}{2} )} = 0 \}^c 
\\ 
& = & 
 A \cup 
 \{ \omega \in \Omega \ : \ 
 \bar{K} (A,\omega ) >0 \} 
 = ((A^c)^\circ)^c = A\cup \partial_{\rm out} A. 
\end{eqnarray*} 
 More refined definitions of inner and outer boundaries 
 are possible, by distinguishing between ``forward'' and ``backward'' 
 neighbors. 
 Note however that defining the norms and boundaries with respect 
 to $K^+$ only, resp. $K^-$ only, leads to 
$\partial_{\rm out} \{ \omega (B ) \leq k \} 
 = \emptyset$ 
 since 
$\vert D^- 1_{\{ \omega (B ) \leq k \}} 
 \vert_{L^p (\sigma )} = 0$, 
 resp. 
 $\partial_{\rm in} \{ \omega (B ) \geq k \} 
 = \emptyset$ 
 since 
$\vert D^+ 1_{\{ \omega (B ) \geq k \}} 
 \vert_{L^p (\sigma )} = 0$, 
 i.e. the isoperimetric constants $h^\pm_p$ 
 defined below have trivial zero value. 
 We have 
$$\pi (\partial_{\rm in} A) 
 = E[ \vert D^+ 1_A\vert_{L^\infty (\sigma+\omega )} ] 
 = \pi (\{ \omega \in A \ : \ 
 \bar{K}(\omega , A^c) >0 \}), 
$$ 
$$\pi (\partial_{\rm out} A) 
 = E[ \vert D^- 1_A\vert_{L^\infty (\sigma+\omega )} ] 
 = \pi (\{ \omega \in A^c \ : \ 
 \bar{K}(\omega , A) >0 \}), 
$$ 
 and 
\begin{eqnarray*} 
\pi (\partial A) 
 & = & E[ \vert D1_A\vert_{L^\infty (\sigma+\omega )} ] 
 = 
 E[ \vert D^+ 1_A\vert_{L^\infty (\sigma+\omega )} ] 
 + E[ \vert D^- 1_A\vert_{L^\infty (\sigma+\omega )} ] 
\\ 
 & = & \pi (\{ \omega \in A \ : \ 
 \bar{K}(\omega , A^c) >0\}) 
 + 
 \pi (\{ \omega \in A^c \ : \ 
 \bar{K}(\omega , A) >0 \}). 
\end{eqnarray*} 
 In discrete settings the surface measure 
 $\pi_s (\partial A)$ 
 of $\partial A$ 
 is not defined via a Minkowski content of the form 
$$\pi_s (\partial A) = \liminf_{r\to 0} 
 \frac{1}{r} (\pi (\{\omega \in \omega \ : \ d(\omega , A)<r\}) - \pi (A)). 
$$ 
 Nevertheless, the surface measure of $\partial_{\rm in} A$, 
 resp. $\partial_{\rm out} A$, can defined by averaging 
 $1_A(\omega ) \bar{K}(\omega , A^c)^{1/2} = 
 \vert D^+ 1_A (\omega ) \vert_{L^2(\frac{\sigma+\omega}{2})}$, 
 resp. $1_{A^c}(\omega ) \bar{K}(\omega, A)^{1/2} 
 = \vert D^- 1_A (\omega ) \vert_{L^2(\frac{\sigma+\omega}{2})}$ 
 with respect to the Poisson measure $\pi (d\omega )$. 
\begin{definition} 
 Let 
$$\pi_s (\partial_{\rm in} A ) 
 = 
 E[ 
 \vert D^+ 1_A (\omega ) \vert_{L^2 (\frac{\sigma+\omega}{2} )} 
 ] 
 = 
 \int_A \bar{K}(\omega , A^c)^{1/2} 
 \pi (d\omega ), 
$$ 
 and 
$$\pi_s (\partial_{\rm out} A ) 
 = 
 E[ 
 \vert D^- 1_A (\omega ) \vert_{L^2 (\frac{\sigma+\omega}{2} )} 
 ] 
 = \int_{A^c} \bar{K}(\omega , A)^{1/2} \pi (d\omega ). 
$$ 
\end{definition} 
 The above quantities 
 represent average numbers of points in $A$, resp. $A^c$, 
 which have a neighbor in $A^c$, resp. $A$, 
 the Poisson measure playing here the role of a uniform measure. 
 The surface measure of $\partial A$ is 
\begin{eqnarray*} 
 \pi_s (\partial A ) & = & 
 \pi_s (\partial_{\rm in} A ) 
 + \pi_s (\partial_{\rm out} A ) 
\\ 
 & = & 
 E\left[ \vert D^+ 1_A \vert_{L^2(\frac{\sigma +\omega}{2})}\right] 
 + E\left[ \vert D^- 1_A \vert_{L^2(\frac{\sigma + \omega}{2} )}\right] 
 = E\left[ \vert D 1_A \vert_{L^2(\frac{\sigma+\omega}{2} )}\right] 
\\ 
& = & 
 \int_A \bar{K}(\omega , A^c)^{1/2} 
 \pi (d\omega ) 
 + \int_{A^c} \bar{K}(\omega , A)^{1/2} 
 \pi (d\omega ). 
\end{eqnarray*} 
 As a consequence of the Margulis-Russo identity Cor.~\ref{margulisrusso} 
 we obtain asymptotic deviation bounds on $\pi_\lambda (A)$ 
 when $A$ is a monotone set. 
\begin{prop} 
 Let $A$ be a monotone subset of $\Omega$, 
 and assume that there exists $\theta > 0$ such that $\pi_\theta (A) = 1/2$. 
 If $A$ is increasing, let 
$$\Delta^- = \inf_{\partial_{\rm out} A} \Vert D^- 1_A\Vert_{L^1 (\sigma )}.$$ 
 We have for $\lambda > \theta$: 
$$ 
\pi_\lambda (A) 
 \leq 
 \Phi  \left(
 \sqrt{2\lambda \Delta^- }-\sqrt{2\theta \Delta^- } 
 \right), 
$$ 
 and for $\lambda < \theta$: 
$$ 
\pi_\lambda (A) 
 \geq 
 \Phi  \left(
 \sqrt{2\lambda \Delta^- }-\sqrt{2\theta \Delta^- } 
 \right). 
$$ 
 If $A$ is decreasing, let 
$$\Delta^+ = \inf_{\partial_{\rm in} A} \Vert D^+ 1_A\Vert_{L^1 (\sigma )},$$ 
 then 
$$ 
\pi_\lambda (A) 
 \leq 
 \Phi  \left(
 \sqrt{2\theta \Delta^+ }-\sqrt{2\lambda \Delta^+ } 
 \right), \quad 
 \lambda > \theta, 
$$ 
 and 
$$ 
\pi_\lambda (A) 
 \geq 
 \Phi  \left(
 \sqrt{2\theta \Delta^+ }-\sqrt{2\lambda \Delta^+ } 
 \right), \quad 
 \lambda < \theta. 
$$ 
\end{prop} 
\begin{Proof} 
 We adapt an argument of \cite{talagrand}, 
 \cite{zemor} to the Poisson case. 
 We have 
\begin{eqnarray*} 
 E_\lambda [\Vert D^- 1_A\Vert_{L^2(\sigma )} ] 
 & = & 
 E_\lambda [1_{\{\Vert D^- 1_A\Vert_{L^\infty (\sigma )}>0\}} 
 \Vert D^- 1_A\Vert_{L^2(\sigma )} ] 
\\ 
 & \leq & 
 \pi_\lambda (\{\Vert D^- 1_A\Vert_{L^\infty (\sigma )}>0\})^{1/2} 
 E_\lambda [\Vert D^- 1_A\Vert_{L^2(\sigma )}^2 ]^{1/2} 
\\ 
 & \leq & 
 \pi_\lambda (\partial_{\rm out} A )^{1/2} 
 E_\lambda [\Vert D^- 1_A\Vert_{L^1(\sigma )} ]^{1/2} 
\\ 
 & \leq & 
 \frac{1}{\sqrt{\Delta^-}} 
 E_\lambda [\Vert D^- 1_A\Vert_{L^1(\sigma )} ]. 
\end{eqnarray*} 
 Let $f(\lambda ) = \pi_\lambda (A)$. 
 Using \eqref{e1} we get 
\begin{eqnarray*} 
 f'(\lambda ) 
 & = & 
 E_\lambda [\Vert D^- 1_A\Vert_{L^1(\sigma )} ]. 
\\ 
 & \geq & 
 \sqrt{\Delta^-} 
 E_\lambda [\Vert D^- 1_A\Vert_{L^2(\sigma )} ]  
\\ 
 & \geq & 
 \sqrt{\frac{\Delta^-}{2\lambda}} I(f(\lambda ))
\\ 
 & = & 
 \frac{ - \sqrt{\Delta^-}}{\sqrt{2\lambda} I''(f(\lambda ))}. 
\end{eqnarray*} 
 Hence for $\lambda > \theta$, 
\begin{eqnarray*} 
\Phi^{-1} (f(\lambda )) 
 & = & \Phi^{-1} (f(\lambda )) - \Phi^{-1} ( f ( \theta ) ) 
\\ 
& = & I'(f(\theta) ) - I'(f(\lambda )) 
\\ 
 & = & 
\int_\lambda^\theta I''(f(t ) ) 
 f'(t ) dt 
\\ 
& \leq & 
 - 
 \int_\lambda^\theta 
 \frac{\sqrt{\Delta^-}}{\sqrt{2 t }} dt 
\\ 
& = & 
 \sqrt{2\Delta^-} 
 (\sqrt{\lambda}-\sqrt{\theta}),  
\end{eqnarray*} 
 and finally 
$$ 
f(\lambda ) 
 \leq 
 \Phi \left(
 \sqrt{2\lambda \Delta^- }-\sqrt{2\theta \Delta^- } 
 \right). 
$$ 
 If $A$ is decreasing and $\lambda > \theta$ we similarly show that 
$$ 
 E_\lambda [\Vert D^+ 1_A\Vert_{L^2(\sigma )} ] 
 \leq 
 \pi_\lambda (\partial_{\rm in} A )^{1/2} 
 E_\lambda [\Vert D^+ 1_A\Vert_{L^1(\sigma )} ]^{1/2} 
 \leq 
 \frac{1}{\sqrt{\Delta^+}} 
 E_\lambda [\Vert D^+ 1_A\Vert_{L^1(\sigma )} ], 
$$ 
\begin{eqnarray*} 
 f'(\lambda ) 
 & = & 
 - E_\lambda [\Vert D^+ 1_A\Vert_{L^1(\sigma )} ]. 
\\ 
 & \leq & 
 \sqrt{\Delta^+} 
 E_\lambda [\Vert D^+ 1_A\Vert_{L^2(\sigma )} ]  
\\ 
 & \leq & 
 \frac{ \sqrt{\Delta^+}}{\sqrt{2\lambda} I''(f(\lambda ))}, 
\end{eqnarray*} 
 and 
$$ 
\Phi^{-1} (f(\lambda )) 
 \leq 
 \int_\lambda^\theta 
 \frac{\sqrt{\Delta^+}}{\sqrt{2 t }} dt 
 = 
 \sqrt{2\Delta^+} 
 (\sqrt{\theta}-\sqrt{\lambda}). 
$$ 
 The case $\lambda < \theta$ is treated in a similar way. 
\end{Proof} 
 When $\lambda<\theta$ and $\Delta^-$ is large, 
 the lower bound is equivalent to 
$$\frac{1}{\sqrt{2\pi} 
 (\sqrt{2\theta \Delta^- }-\sqrt{2\lambda \Delta^- })} 
 e^{-(\sqrt{2\lambda \Delta^- }-\sqrt{2\theta \Delta^- })^2/2} 
. 
$$ 
 As an example, for the increasing set $\{\omega (B) \geq n \}$ we have 
$$\partial_{\rm out} \{\omega (B) \geq n \} = \{\omega (B) = n-1 \},$$ 
 and 
$$D_x 1_{\{\omega (B) \geq n \}} = 
 - D^-_x 1_{\{\omega (B) \geq n \}} = 
 - 1_B (x) 
 1_{\{\omega (B) = n-1 \}}, 
$$ 
 hence 
$$\Vert D 1_{\{\omega (B) \geq n \}} \Vert_{L^1(\sigma ) } 
 = \sigma (B) 1_{\{\omega (B) = n-1 \}} 
 = \sigma (B) 1_{\partial_{\rm out} \{ \omega (B) \geq n \}},$$ 
 and $\Delta^- = \sigma (B)$. 
 For the decreasing set $\{\omega (B) \leq n \}$ we have 
$$\partial_{\rm in} \{\omega (B) \leq n \} = \{\omega (B) = n \},$$ 
 and 
$$D_x 1_{\{\omega (B) \leq n \}} = 
 - D^-_x 1_{\{\omega (B) \leq n \}} = 
 - 1_B (x) 
 1_{\{\omega (B) = n \}}, 
$$ 
 hence 
$$\Vert D 1_{\{\omega (B) \leq n \}} \Vert_{L^1(\sigma ) } 
 = \sigma (B) 1_{\{\omega (B) = n \}} 
 = \sigma (B) 1_{\partial_{\rm in} \{ \omega (B) = n \}},$$ 
 and $\Delta^+ = \sigma (B)$. 
\section{Co-area formulas} 
\label{s4} 
 For $p=\infty$ the next Lemma shows that 
$$E[\vert D^+ F \vert_{L^\infty (\sigma+\omega )}] 
 = \int_{-\infty}^{+\infty} 
 \pi ( 
 \partial_{\rm in} \{ F>t\} ) 
 dt,
$$ 
$$E[\vert D^- F \vert_{L^\infty (\sigma+\omega )} ] 
 = \int_{-\infty}^{+\infty} 
 \pi ( 
 \partial_{\rm out} \{ F>t\} 
 ) 
 dt. 
$$ 

\begin{lemma} 
\label{co-area} 
 We have 
$$E[\vert D^\pm F \vert_{L^\infty (\sigma+\omega )}] 
 = \int_{-\infty}^{+\infty} 
 E[ \vert D^\pm 1_{\{F>t\}} \vert_{L^\infty (\sigma+\omega )}] 
 dt 
, 
$$ 
 and 
$$E\left[\vert D^+ F\vert_{L^\infty (\sigma+\omega )}\right] 
 + E\left[\vert D^- F\vert_{L^\infty (\sigma+\omega )}\right] 
 = \int_{-\infty}^\infty 
 E\left[\vert D 1_{\{F>t\}}\vert_{L^\infty (\sigma+\omega )} 
  \right] 
 dt. 
$$ 
\end{lemma} 
\begin{Proof} 
 The notations 
 $\esssup_{\tilde{\omega } \in {\cal N}_\omega}$ 
 and 
 $\essinf_{\tilde{\omega } \in {\cal N}_\omega}$ 
 denote respectively 
 $\esssup_{\bar{K}(\omega, d \tilde{\omega })}$ 
 and 
 $\essinf_{\bar{K}(\omega, d \tilde{\omega })}$. 
 We have 
 $$
 \vert D^+ F(\omega ) \vert_{L^\infty (\sigma+\omega )} 
 = \esssup_{\tilde{\omega } \in {\cal N}_\omega} 
 ( F(\omega ) - F(\tilde{\omega} ))^+ 
 = F(\omega ) 
 - \essinf_{\tilde{\omega } \in {\cal N}_\omega} 
 F(\tilde{\omega}), 
$$ 
 hence 
\begin{eqnarray*} 
\lefteqn{ 
E[\vert D^+ F \vert_{L^\infty (\sigma+\omega )}] 
 = E[F] - E[
 \essinf_{\tilde{\omega } \in {\cal N}_\omega}  
 F(\tilde{\omega})] 
} 
\\ 
& = & 
\int_{-\infty}^{+\infty} \pi (\{F>t\}) dt 
- \int_{-\infty}^{+\infty} 
  \pi ( 
 \essinf_{\tilde{\omega } \in {\cal N}_\omega}  
 F(\tilde{\omega})>t) 
 dt 
\\ 
& = & 
\int_{-\infty}^{+\infty} \pi (\{F>t\}) dt 
- \int_{-\infty}^{+\infty} 
  \pi (\{ 
 \essinf_{\tilde{\omega } \in {\cal N}_\omega} 
  F(\tilde{\omega}) >t \ \mbox{and} 
 \ F(\omega ) >t\} ) dt 
\\ 
& = & 
\int_{-\infty}^{+\infty} \pi (\{F (\omega ) >t 
 \ \mbox{and} 
 \ (\sigma+\omega )(\{ x\in X \ : \ F(\omega \pm \delta_x ) \leq t\})>0 \}) dt  
\\ 
& = & 
\int_{-\infty}^{+\infty} \pi (\{ 
 \omega \in \Omega \ : \ 
 (\sigma+\omega ) (\{ x \in X \ : 
 F(\omega )>t \ \mbox{and} 
 \  F(\omega \pm \delta_x ) \leq t\})>0 \}) dt  
\\ 
& = & 
\int_{-\infty}^{+\infty} \pi (\{ 
 \omega \in \Omega \ : \ 
 (\sigma+\omega ) (\{ x \in X \ : 
 1_{\{F(\omega )>t\}}-1_{\{F(\omega \pm \delta_x ) > \}}=1\} 
 )>0 \}) dt  
\\ 
& = & 
\int_{-\infty}^{+\infty} \pi (\{ 
 \omega \in \Omega \ : \ 
 \vert D^+ 1_{\{F>t\}} \vert_{L^\infty (\sigma+\omega )}=1 
 \}) dt  
\\ 
& = & 
\int_{-\infty}^{+\infty} E[\vert D^+ 1_{\{F>t\}} \vert_{L^\infty (\sigma+\omega )} ] 
 dt. 
\end{eqnarray*} 
 The proof for $D^-$ is similar. 
 Finally we have, since $D_x^+FD_x^-F=0$: 
\begin{eqnarray*} 
\lefteqn{ 
 E[\vert D^+ F\vert_{L^\infty (\sigma+\omega )}] 
 + E[\vert D^- F\vert_{L^\infty (\sigma+\omega )}] 
} 
\\ 
 & = & 
 \int_{-\infty}^\infty 
 E[\vert D^+1_{\{F>t\}}\vert_{L^\infty (\sigma+\omega )}] 
 dt 
 + 
 \int_{-\infty }^\infty 
 E[\vert D^- 1_{\{F>t\}}\vert_{L^\infty (\sigma+\omega )}] 
 dt 
\\ 
 & = & 
 \int_{-\infty}^\infty 
 E[\vert D^+1_{\{F>t\}}\vert_{L^\infty (\sigma+\omega )} 
  + \vert D^- 1_{\{F>t\}}\vert_{L^\infty (\sigma+\omega )}] 
 dt 
\\ 
 & = & 
 \int_{-\infty}^\infty 
 E[\vert D 1_{\{F>t\}}\vert_{L^\infty (\sigma+\omega )} 
  ] 
 dt. 
\end{eqnarray*} 
\end{Proof} 
 The next Lemma states a co-area formula in $L^1$. 
\begin{lemma} 
\label{co-area2} 
 We have 
$$E[\vert D^\pm F \vert_{L^1 (\sigma )}] 
 = \int_{-\infty}^{+\infty} 
 E[ \vert D^\pm 1_{\{F>t\}} \vert_{L^1 (\sigma )}] 
 dt, 
$$ 
$$E[\vert D^\pm F \vert_{L^1 (\omega )}] 
 = \int_{-\infty}^{+\infty} 
 E[ \vert D^\pm 1_{\{F>t\}} \vert_{L^1 (\omega )}] 
 dt. 
$$ 
\end{lemma} 
\begin{Proof} 
 We have for all $a,b\in \real$: 
$$(b-a)^\pm = \int_{-\infty}^\infty 
 (1_{\{a>t\}} 
 - 1_{\{b>t\}} )^\pm dt,$$ 
 hence 
$$ 
D^\pm_x F 
 = \int_{-\infty}^\infty 
 D^\pm_x 1_{\{F >t\}} dt.$$ 
\end{Proof} 
 As a consequence we have 
$$E\left[ 
 \vert D^\pm F \vert_{L^1 (\frac{\sigma +\omega}{2} )} 
 \right] 
 = \int_{-\infty}^{+\infty} 
 E\left[ \vert D^\pm 1_{\{F>t\}} \vert_{L^1 (\frac{\sigma + \omega}{2} )}\right] 
 dt, 
$$ 
 and 
$$E\left[ 
 \vert D F \vert_{L^1 (\frac{\sigma +\omega}{2} )} 
 \right] 
 = \int_{-\infty}^{+\infty} 
 E\left[ \vert D 1_{\{F>t\}} \vert_{L^1 (\frac{\sigma + \omega}{2} )}\right] 
 dt. 
$$ 
\begin{prop} 
 We have 
$$E[\Gamma^\pm (F,F)] 
 = \int_{-\infty}^{+\infty} \int_{-\infty}^{+\infty} 
 E[ \Gamma^\pm ( 1_{\{F>t\}} ,  1_{\{F>s\}})] ds dt. 
$$ 
\end{prop} 
\begin{Proof} 
 Again we use the identity 
$$D^\pm_x F = \int_{-\infty}^{+\infty} 
 D^\pm_x 1_{\{F>t\}} dt.
$$ 
\end{Proof} 
 We close this section with an application of co-area formulas 
 to an equivalence result on functional inequalities. 
 Let ${\cal G}$ be a non-empty set of functions on $\Omega$, and let 
\begin{equation} 
\label{eq} 
{\cal L} (F) = \sup_{G_1,G_2 \in {\cal G}} 
 E[F^+G_1+F^-G_2].
\end{equation} 
 Several functionals have the representation \eqref{eq}, 
 for example the entropy 
$${\cal L}(F) 
 = \Ent \vert F \vert = E[\vert F \vert 
 \log \vert F \vert ] 
 - E[\vert F\vert ]\log E[\vert F \vert ] 
 = \sup_{E[e^G]\leq 1} E[\vert F \vert G], 
$$ 
 the variance 
$${\cal L}(F) 
 = E[(F-E[F])^2] = \var (F) 
 = \inf_{a\in \real} E[( F-a)^2 ] 
, 
$$ 
 and 
$${\cal L}(F) 
 = E[\vert F- m(F)\vert ] 
 = \inf_{a\in \real} E[\vert F-a\vert ] 
, 
$$ 
 where $m(F)$ is by definition a median of $F$. 
 The co-area formula implies the following 
 equivalence, as in \cite{houdremixed}, \cite{rothaus}. 
 The norm $\vert \cdot \vert_p$ denotes 
 either $\vert \cdot \vert_{L^1 (\sigma )}$ 
 or $\vert \cdot \vert_{L^1 (\omega )}$ 
 when $p=1$, and 
 $\vert \cdot \vert_{L^\infty (\sigma+\omega )}$ 
 when $p=\infty$. 
\begin{theorem}
\label{iso} 
 Let $c\geq 0$. The following are equivalent: 
\begin{description} 
\item{(i)} $c {\cal L} (F) \leq E[ 
 \vert D^\pm F \vert_p]$, for all $F:\Omega \rightarrow \real$, 
\item{(ii)} $c {\cal L} (1_A) \leq E[ 
 \vert D^\pm 1_A \vert_p]$ 
 and 
 $c {\cal L} (-1_A) \leq E[ 
 \vert D^\pm (-1_A ) \vert_p]$, 
 for all $A\in {\cal F}$, 
\end{description} 
 with $p=1,\infty$. 
\end{theorem} 
\begin{Proof} 
 We follow the proof of \cite{houdremixed}. 
 In order to show $(ii) \Rightarrow (i)$ 
 we note that for all $G_1,G_2\in {\cal G}$, 
\begin{eqnarray*} 
 E[ \vert D^\pm F \vert_p] & = & 
 \int_0^\infty 
 E[ \vert D^\pm 1_{\{F>t\}} \vert_p] dt 
 + 
 \int_{-\infty}^0 
 E[ \vert D^\pm 1_{\{F>t\}} \vert_p] dt 
\\ 
 & \geq & 
 c \int_0^\infty 
 E[ G_1 1_{\{F>t\}} ] dt 
 + \int_{-\infty}^0 
 E[ \vert D^\pm (-1_{\{F\leq t\}}) \vert_p] dt 
\\ 
 & \geq & 
 c E[ G_1 F^+ ] + c \int_{-\infty}^0 
 E[ 1_{\{F\leq t\}} G_2] dt 
\\ 
 & = & 
 c E[ G_1 F^+ ] + c 
 E[ F^- G_2], 
\end{eqnarray*} 
 hence 
$$ 
 E[ \vert D^\pm F \vert_p] 
 \geq 
 c \sup_{G_1,G_2\in {\cal G}} 
 ( 
 E[ G_1 F^+ ] + E[ F^- G_2] 
 ) 
 \geq c {\cal L}(F). 
$$ 
\end{Proof} 
\section{Some explicit computations} 
\label{s7} 
 In this section we define the main isoperimetric 
 constants and establish some bounds on these constants. 
\begin{definition} 
 Let for $p\in [1,\infty ]$: 
$$h^\pm_p 
 = \inf_{0 < \pi (A ) < \frac{1}{2}} 
 \frac{E\left[ 
 \vert D^\pm 1_A\vert_{L^p(\frac{\sigma+\omega}{2} )}\right] 
 }{\pi (A)}, 
 \quad 
 \quad  
 h_p = \inf_{0 < \pi (A ) < \frac{1}{2}} 
 \frac{E\left[ 
 \vert D 1_A\vert_{L^p(\frac{\sigma+\omega}{2} )}\right] 
 }{\pi (A)}. 
$$ 
\end{definition} 
 We have 
$$ 
 h^+_1  
 = h^-_1  
 = \inf_{0 < \pi (A ) < \frac{1}{2}} 
 \frac{1}{\pi (A)} 
 \int_A \bar{K}(\omega , A^c)\pi (d\omega ) 
 = \inf_{0 < \pi (A ) < \frac{1}{2}} 
 \frac{1}{\pi (A)} 
 \int_{A^c} \bar{K}(\omega , A ) \pi (d\omega ), 
$$
$$ 
h^+_2 
 = \inf_{0 < \pi (A ) < \frac{1}{2}} 
 \frac{\pi_s (\partial_{\rm in} A)}{\pi (A)}, 
\quad 
h^-_2 
 = \inf_{0 < \pi (A ) < \frac{1}{2}} 
 \frac{\pi_s (\partial_{\rm out} A)}{\pi (A)}, 
\quad 
h_2 = \inf_{0 < \pi (A ) < \frac{1}{2}} 
 \frac{\pi_s (\partial A)}{\pi (A)} 
$$ 
 and 
$$ 
h^+_\infty 
 = \inf_{0 < \pi (A ) < \frac{1}{2}} 
 \frac{\pi (\partial_{\rm in} A)}{\pi (A)}, 
\quad 
h^-_\infty 
 = \inf_{0 < \pi (A ) < \frac{1}{2}} 
 \frac{\pi (\partial_{\rm out} A)}{\pi (A)}, 
\quad 
h_\infty = \inf_{0 < \pi (A ) < \frac{1}{2}} 
 \frac{\pi (\partial A)}{\pi (A)}. 
$$ 
 The following is a functional version of $h^\pm_p$. 
\begin{definition} 
 Let for $p\in [1,\infty ]$: 
$$\tilde{h}^\pm_p 
 = \inf_{0 < \pi (A ) < 1} 
 \frac{E\left[ 
 \vert D^\pm 1_A\vert_{L^p(\frac{\sigma+\omega}{2} )}\right] 
 }{\pi (A)\pi (A^c)}, 
 \quad 
 \quad  
 \tilde{h}_p = \inf_{0 < \pi (A ) < 1} 
 \frac{E\left[ 
 \vert D 1_A\vert_{L^p(\frac{\sigma+\omega}{2} )}\right] 
 }{\pi (A)\pi (A^c)}. 
$$ 
\end{definition} 
 Note that in the definition of the isoperimetric constants 
 we need to integrate with respect to $\omega + \sigma$, 
 otherwise integrating with respect to $\omega$ or 
 $\sigma$ only would lead to vanishing isoperimetric constants, 
 since 
$$ \{ \vert D^+ 1_{\{ \omega (B ) \leq k \}} 
 \vert_{L^p (\omega )}>0\} 
 = \emptyset,$$ 
and 
$$ \{ \vert D^- 1_{\{ \omega (B ) \geq k \}} 
 \vert_{L^p (\sigma )}>0\} = \emptyset. 
$$ 
 The next proposition follows the presentation of \cite{stoyanov}. 
\begin{prop} 
 We have 
\begin{description} 
\item{a)} $h_1=2h^+_1=2h^-_1$, 
\item{b)} $\tilde{h}^+_p=\tilde{h}^-_p$, $p=1,\infty$, 
\item{c)} $\min (h^+_p,h^-_p) 
 < \tilde{h}^\pm_p < 2\min (h_p^+,h_p^-) 
 \leq h^+_p+h^-_p \leq h_p < \tilde{h}_p 
 < 2h_p$, \ $p\in [1,+\infty ]$. 
\end{description} 
\end{prop} 
\begin{Proof} 
 For the first statement we use Lemma~\ref{l1}, which implies 
$$E\left[\vert D^+ F \vert_{L^p(\frac{\sigma+\omega}{2} )}^p \right]  
 = E\left[\vert D^- F \vert_{L^p(\frac{\sigma+\omega}{2} )}^p \right] 
 = \frac{1}{2} 
 E\left[\vert D F \vert_{L^p(\frac{\sigma+\omega}{2} )}^p \right]. 
$$ 
 The second statement follows from Remark~\ref{r1}. 
 The last statement follows from the inequalities, if $0<\pi (A) < 1/2$: 
\begin{eqnarray*} 
 h^+_p & \leq & 
 \frac{E\left[\vert D^+1_A\vert_{L^p(\frac{\sigma+\omega}{2} )}\right]}{\pi (A)} 
 \leq 
 \frac{E\left[\vert D^+1_A\vert_{L^p(\frac{\sigma+\omega}{2} )}\right]}{\pi (A) 
 \pi (A^c)} 
\\ 
& = & 
 \frac{E\left[\vert D^\pm 1_A\vert_{L^p(\frac{\sigma+\omega}{2} )}\right]}{\pi (A) 
 \pi (A^c)} 
 \leq 
 2 \frac{E\left[\vert D^\pm 1_A\vert_{L^p(\frac{\sigma+\omega}{2} )}\right]}{\pi (A)} 
, 
\end{eqnarray*} 
 and similarly if $1/2\leq \pi (A) < 1$: 
\begin{eqnarray*} 
h^-_p & \leq & 
 \frac{E\left[\vert D^- 1_{A^c}\vert_{L^p(\frac{\sigma+\omega}{2} )}\right]}{\pi (A^c )} 
 \leq 
 \frac{E\left[\vert D^- 1_{A^c}\vert_{L^p(\frac{\sigma+\omega}{2} )}\right]}{\pi (A) 
 \pi (A^c)} 
\\ 
& = & 
 \frac{E\left[\vert D^\pm 1_{A^c} \vert_{L^p(\frac{\sigma+\omega}{2} )}\right]}{\pi (A) 
 \pi (A^c)} 
 \leq 
 2 \frac{E\left[\vert D^\pm 1_{A^c} \vert_{L^p(\frac{\sigma+\omega}{2} )}\right]}{\pi (A^c)} 
. 
\end{eqnarray*} 
\end{Proof} 
\begin{definition} 
 Let for $p\in [1,\infty ]$: 
$$k^\pm_p = 
 \inf_{F\not= C} \frac{E\left[\vert D^\pm F\vert_{L^p(\frac{\sigma+\omega}{2} )} \right]}{E[( F-m(F))^\pm  ]}, 
 \quad \quad 
k_p = 
 \inf_{F\not= C} \frac{E\left[\vert D F\vert_{L^p(\frac{\sigma+\omega}{2} )} \right]}{E[\vert F-m(F)\vert ]}. 
$$ 
\end{definition} 
\begin{prop} 
 We have 
 $h_1^\pm=k_1^\pm$, $h_\infty^\pm=k_\infty^\pm$, 
 $h_1=k_1$, and $k_\infty \leq h_\infty \leq 2k_\infty$. 
\end{prop} 
\begin{Proof} 
 First of all we note that 
 since $m(1_A) = 0$ if $\pi (A)\leq 1/2$, we have 
$$k^\pm_p \pi (A) 
 = k^\pm_p E[(1_A-m(1_A))^\pm ] 
 \leq 
 E\left[\vert D^\pm 1_A\vert_{L^p (\frac{\sigma+\omega}{2} )}\right] 
, 
$$ 
 hence $h^\pm_p\geq k^\pm_p$, $p=1,\infty$, 
 and similarly 
$$k_p \pi (A) 
 = k_p E[\vert 1_A-m(1_A)\vert ] 
 \leq 
 E\left[\vert D 1_A\vert_{L^p (\frac{\sigma+\omega}{2} )}\right] 
, 
$$ 
 hence $h_p\geq k_p$, $p=1,\infty$, 
 From the co-area formulas Lemmas~\ref{co-area} and \ref{co-area2} 
we have for $p=1,\infty$, 
since $\pi (F>m(F))\leq 1/2$: 
\begin{eqnarray*} 
 E\left[\vert D^+F\vert_{L^p (\frac{\sigma+\omega}{2} )}\right] 
 & = & \int_{-\infty}^\infty 
 E\left[\vert D^+1_{\{F>t\}}\vert_{L^p(\frac{\sigma+\omega}{2} )}\right] 
 dt 
\\ 
 & \geq & 
 h^+_p  
 \int_{-\infty}^\infty 
 \pi ( \{F>t\} ) 
 dt 
\\ 
 & \geq & 
 h^+_p  
 \int_{m(F)}^\infty 
 \pi ( \{F>t\} ) 
 dt 
\\ 
 & = & 
 h^+_p  
 \int_0^\infty 
 \pi ( \{F-m(F)>t\} ) 
 dt 
\\ 
 & \geq & h^+_p E[(F-m(F))^+]. 
\end{eqnarray*} 
 Hence $k^+_p \geq h^+_p$. 
 Similarly we obtain 
$$E\left[\vert D^-F\vert_{L^p (\frac{\sigma+\omega}{2} )}\right] 
 = E\left[\vert D^+(-F)\vert_{L^p (\frac{\sigma+\omega}{2} )}\right] 
 \geq h^-_p E[(-F-m(-F))^+] 
 = h^-_p E[(F-m(F))^-],
$$ 
 hence $k^-_p \geq h^-_p$, and 
\begin{eqnarray*} 
 E\left[\vert D F\vert_{L^1 (\frac{\sigma+\omega}{2} )}\right] 
 & = & \int_{-\infty}^{m(F)} 
 E\left[\vert D1_{\{F>t\}}\vert_{L^1(\frac{\sigma+\omega}{2} )}\right] 
 dt 
 + 
 \int_{m(F)}^\infty 
 E\left[\vert D1_{\{F>t\}}\vert_{L^1(\frac{\sigma+\omega}{2} )}\right] 
 dt 
\\ 
 & = & \int_0^\infty 
 E\left[\vert D1_{\{-F+m(F)>t\}}\vert_{L^1(\frac{\sigma+\omega}{2} )}\right] dt 
 + 
 \int_0^\infty 
 E\left[\vert D1_{\{F-m(F)>t\}}\vert_{L^1(\frac{\sigma+\omega}{2} )}\right] 
 dt 
\\ 
 & \geq & 
 h_1 
 \int_{-\infty}^0 
 \pi ( \{-F+m(F)>t\} ) 
 dt 
 + 
 h_1 
 \int_0^\infty 
 \pi ( \{F-m(F)>t\} ) 
 dt 
\\ 
 & \geq & 
 h_1 
 E[(F-m(F))^-] + 
 h_1 
 E[(F-m(F))^+] 
\\ 
 & = & h_1 E[\vert F-m(F)\vert], 
\end{eqnarray*} 
 hence $k_1\geq h_1$. 
 From Lemma~\ref{co-area} we also have 
\begin{eqnarray*} 
 2 E[\vert D F\vert_{L^\infty (\sigma+\omega )}] 
 & \geq & 
 E[\vert D^+ F\vert_{L^\infty (\sigma+\omega )}] 
 + E[\vert D^- F\vert_{L^\infty (\sigma+\omega )}] 
\\ 
 & = & 
 \int_{-\infty}^\infty 
 E[\vert D 1_{\{F>t\}}\vert_{L^\infty (\sigma+\omega )} 
  ] 
 dt 
\\ 
 & = & 
 \int_{-\infty}^{m(F)} 
 E[\vert D1_{\{F>t\}}\vert_{L^\infty (\sigma+\omega )}] 
 dt 
 + 
 \int_{m(F)}^\infty 
 E[\vert D1_{\{F>t\}}\vert_{L^\infty (\sigma+\omega )}] 
 dt 
\\ 
 & = & \int_0^\infty 
 E[\vert D1_{\{-F+m(F)>t\}}\vert_{L^\infty (\sigma+\omega )}] dt 
 + 
 \int_0^\infty 
 E[\vert D1_{\{F-m(F)>t\}}\vert_{L^\infty (\sigma+\omega )}] 
 dt 
\\ 
 & \geq & 
 h_\infty 
 \int_{-\infty}^0 
 \pi ( \{-F+m(F)>t\} ) 
 dt 
 + 
 h_\infty 
 \int_0^\infty 
 \pi ( \{F-m(F)>t\} ) 
 dt 
\\ 
 & \geq & 
 h_\infty 
 E[(F-m(F))^-] + 
 h_\infty  
 E[(F-m(F))^+] 
\\ 
 & \geq & 
 h_\infty 
 E[\vert F-m(F)\vert], 
\end{eqnarray*} 
 hence $2k_\infty \geq h_\infty$. 
\end{Proof} 
\begin{remark} 
\label{c11} 
 The above proof also implies, if $F\geq 0$ 
 and $\pi (F>0)\leq 1/2$: 
$$h^+_p 
 E[F] 
 \leq E\left[\vert D^+ F\vert_{L^p(\frac{\sigma+\omega}{2} )}\right] 
 , 
$$ 
 and 
$$ 
 h^-_p E[F] 
 \leq E\left[\vert D^- F\vert_{L^p (\frac{\sigma+\omega}{2} )}\right]. 
$$ 
\end{remark} 
 The following is the definition of the Poincar\'e constants. 
\begin{definition} 
 Let for $p\in [1,\infty ]$: 
$$\lambda^\pm_p 
 = 
 \inf_{F\not= C} \frac{E\left[\vert D^\pm F\vert_{L^p(\frac{\sigma+\omega}{2} )}^2 \right]}{\var (F)}, 
 \quad 
 \quad 
\lambda_p 
 = 
 \inf_{F\not= C} \frac{E\left[\vert D F\vert_{L^p(\frac{\sigma+\omega}{2} )}^2 \right]}{\var (F)}. 
$$ 
\end{definition} 
 Remark that $\lambda^+_p=\lambda^-_p$, 
 $p\in [1,\infty ]$, since $D^+_xF = D^-_x(-F)$, 
 and $\tilde{h}^+_1 \geq \lambda^+_2$. 
 We have 
$$ 
E\left[\vert DF\vert_{L^2(\frac{\sigma+\omega}{2} )}^2 \right] 
= \frac{1}{2} 
 E\left[\vert DF\vert_{L^2(\sigma)}^2 \right] 
= \frac{1}{2} 
 E\left[\vert DF\vert_{L^2(\omega)}^2 \right], 
$$ 
 hence 
$$\lambda_2 = 2 \inf_{F\not= C}
 \frac{{\cal E}(F,F)}{\var (F)}. 
$$ 
 Th.~\ref{iso} shows that 
$$\lambda^\pm_\infty = 
 \inf_{\pi (A ) >0 } 
 \frac{E[\vert D^\pm 1_A\vert_{L^\infty (\sigma+\omega )}]}{\var 1_A}. 
$$ 

 \begin{definition} 
 Let for $p\in [1,\infty ]$: 
$$\tilde{k}^\pm_p = 
 \inf_{F\not= C} \frac{E[\vert D^\pm F\vert_{L^p(\frac{\sigma+\omega}{2} )} ]}{E[ ( F-E[F])^\pm  ]}, 
 \quad \quad 
 \tilde{k}_p = 
 \inf_{F\not= C} \frac{E[\vert D F\vert_{L^p(\frac{\sigma+\omega}{2} )} ]}{E[\vert F-E[F] \vert ]}. 
$$ 
\end{definition} 
\begin{prop} 
\label{6.6} 
 We have 
$$k^+_\infty = h^+_\infty 
 \leq 2 \sqrt{\lambda_\infty} 
 = \frac{2}{\sqrt{\sigma (X)}} \quad 
 \mbox{and} 
 \quad 
 k^-_\infty = h^-_\infty \leq \frac{2}{\sigma (X)}\left( 
 1 + \sqrt{2\sigma (X)}\right). 
$$ 
\end{prop} 
\begin{Proof} 
 Note that if $F\geq 0$, 
\begin{eqnarray*} 
\vert D^+ F^2 (\omega ) 
 \vert_{L^\infty (\sigma+\omega )} 
 & = & 
 \esssup_{\tilde{\omega } \in {\cal N}_\omega} 
 (F^2 (\omega ) - F^2 (\tilde{\omega} )) 
\\ 
 & = & 
 \esssup_{\tilde{\omega } \in {\cal N}_\omega} 
 (F^2 (\omega ) - F^2 (\tilde{\omega} )) 
 1_{\{F (\omega ) \geq F (\tilde{\omega} )\}} 
\\ 
 & = & 
 \esssup_{\tilde{\omega } \in {\cal N}_\omega}  
 (F (\omega ) - F (\tilde{\omega} )) 
 (F (\omega ) + F (\tilde{\omega} )) 
 1_{\{F (\omega ) \geq F (\tilde{\omega} )\}} 
\\ 
 & \leq & 
 2 
 \esssup_{\tilde{\omega } \in {\cal N}_\omega} 
 (F (\omega ) - F (\tilde{\omega} )) 
 F (\omega ) 
\\ 
 & = & 
 2 \vert D^+ F\vert_{L^\infty (\sigma+\omega )} F ( \omega ). 
\end{eqnarray*} 
 If $\pi (\{ F>0\}) \leq 1/2$, then by Remark~\ref{c11} applied 
 to $F^2$, 
\begin{eqnarray*} 
 (h^+_\infty)^2 E[F^2]^2 
 & \leq & 
 E[\vert D^+ F^2 \vert_{L^\infty (\sigma )}]^2 
\\ 
 & \leq & 
 4 E[F \vert D^+ F\vert_{L^\infty (\sigma )}]^2 
\\ 
 & \leq & 
 4 E[\vert D^+ F\vert_{L^\infty (\sigma+\omega )}^2] 
 E[F^2],  
\end{eqnarray*} 
 hence 
$$ 
 \frac{(h^+_\infty)^2}{4} 
 E[F^2] 
 \leq E[\vert D^+ F \vert_{L^\infty (\sigma )}^2]. 
$$ 
 In the general case we may assume 
 that $m(F) = 0$, i.e. 
$$\pi (\{ F>0\}) \leq 1/2, \quad \mbox{and} 
 \quad \pi (\{ F<0\}) \leq 1/2. 
$$ 
 We have 
$$\pi (\{ F^+>0\}) \leq 1/2, \quad \mbox{and} 
 \quad \pi (\{ F^- <0\}) \leq 1/2, 
$$ 
 hence 
$$\frac{(h^+_\infty)^2}{4} 
 E[(F^+)^2] 
 \leq E[\vert D^+ F^+ \vert_{L^\infty (\sigma )}^2], 
$$ 
 and 
$$\frac{(h^+_\infty)^2}{4} 
 E[(F^-)^2] 
 \leq E[ \vert D^+ F^- \vert_{L^\infty (\sigma )}^2]. 
$$ 
 We have 
\begin{eqnarray*} 
\vert D^+ F^+ (\omega ) \vert_{L^\infty (\sigma )} 
 & = & 
\esssup_{\tilde{\omega } \in {\cal N}_\omega} 
 (F^+ (\omega ) - F^+ (\tilde{\omega} )) 
\\ 
 & = & 
\esssup_{\tilde{\omega } \in {\cal N}_\omega} 
 \vert 
 F (\omega ) - F (\tilde{\omega} ) \vert 
 1_{\{F (\omega ) > 0\}}, 
\end{eqnarray*} 
 and 
\begin{eqnarray*} 
\vert D^+ F^- (\omega ) \vert_{L^\infty (\sigma )} 
 & = & 
\esssup_{\tilde{\omega } \in {\cal N}_\omega} 
 (F^- (\omega ) - F^- (\tilde{\omega} )) 
\\ 
 & = & 
\esssup_{\tilde{\omega } \in {\cal N}_\omega} 
 \vert F (\omega ) - F (\tilde{\omega} )\vert 
 1_{\{F (\omega ) < 0\}}. 
\end{eqnarray*} 
 Hence 
\begin{eqnarray*} 
 \frac{(h^+_\infty)^2}{4} 
 \var F & \leq &  \frac{(h^+_\infty)^2}{4} 
 E[F^2] 
 = \frac{(h^+_\infty)^2}{4} 
 E[F^21_{\{F>0\}}] 
 + \frac{(h^+_\infty)^2}{4} 
 E[F^21_{\{F<0\}}] 
\\ 
 & \leq & 
 E[1_{\{F>0\}} \vert D F\vert_{L^\infty (\sigma )}^2] 
 + E[1_{\{F<0\}} \vert D F\vert_{L^\infty (\sigma )}^2], 
\end{eqnarray*} 
 from which $\lambda_\infty \geq (h^+_\infty)^2/4$. 
 The second statement has a similar proof: 
\begin{eqnarray*} 
\vert D^- F^2 \vert_{L^\infty (\sigma )} 
 & = & 
\esssup_{\tilde{\omega } \in {\cal N}_\omega} 
 (F^2 (\tilde{\omega} ) - F^2 (\omega )) 
\\ 
 & = & 
\esssup_{\tilde{\omega } \in {\cal N}_\omega} 
 (F (\tilde{\omega} ) - F (\omega )) 
 (F (\tilde{\omega} ) + F (\omega )) 
 1_{\{F (\tilde{\omega} ) \geq F (\omega )\}} 
\\ 
 & = & 
\esssup_{\tilde{\omega } \in {\cal N}_\omega} 
 (F (\tilde{\omega} ) - F (\omega ))^2 
 + 2(F (\tilde{\omega} ) - F (\omega )) 
 F (\omega ) 1_{\{F (\tilde{\omega} ) \geq F (\omega )\}} 
\\ 
 & \leq & 
\esssup_{\tilde{\omega } \in {\cal N}_\omega} 
 (F (\tilde{\omega} ) - F (\omega ))^2 
 + 2(F (\tilde{\omega} ) - F (\omega )) 
 F (\omega ) 1_{\{F (\tilde{\omega} ) \geq F (\omega )\}}. 
\end{eqnarray*} 
 By Remark~\ref{c11}, 
\begin{eqnarray*}  
 h^-_\infty E[F^2] 
 & \leq & 
 E\left[ 
 \esssup_{\tilde{\omega } \in {\cal N}_\omega} 
 (F (\tilde{\omega} ) - F (\omega ))^2 
 + 2(F (\tilde{\omega} ) - F (\omega )) 
 F (\omega ) 1_{\{F (\tilde{\omega} ) \geq F (\omega )\}}\right] 
\\ 
 & \leq & 
 E\left[ 
 \esssup_{\tilde{\omega } \in {\cal N}_\omega} 
 (F (\tilde{\omega} ) - F (\omega ))^2 
 \right]  
\\ 
 & & + 
 2 
 E[F^2]^{1/2} 
 E\left[ 
 \esssup_{\tilde{\omega } \in {\cal N}_\omega} 
 (F (\tilde{\omega} ) - F (\omega ))^2 
 1_{\{F (\tilde{\omega} ) \geq F (\omega )\}}\right]^{1/2},  
\end{eqnarray*} 
 hence 
\begin{equation} 
\label{ll2} 
 (\sqrt{1+h^-_\infty}-1)^2 
 E[F^2] 
 \leq 
 E\left[ 
 \esssup_{\tilde{\omega } \in {\cal N}_\omega} 
 (F (\tilde{\omega} ) - F (\omega ))^2 
 1_{\{F(\tilde{\omega } ) \geq F(\omega )\}}  
 \right]. 
\end{equation} 
 In the general case, if $0$ is a median of $F$ 
 we have, applying \eqref{ll2} to $F^+$ and $F^-$: 
\begin{eqnarray*} 
E[F^2] & = & 
E[(F^+)^2] 
+ E[(F^-)^2] 
\\ 
& \leq & 
 E\left[ 
 \esssup_{\tilde{\omega } \in {\cal N}_\omega} 
 (F^+ (\tilde{\omega} ) - F^+ (\omega ))^2 
 1_{\{F^+ (\tilde{\omega } ) \geq F^+ (\omega )\}} 
 \right] 
\\ 
 & & + E\left[ 
 \esssup_{\tilde{\omega } \in {\cal N}_\omega} 
 (F^- (\tilde{\omega} ) - F^- (\omega ))^2 
 1_{\{F^- (\tilde{\omega } ) \geq F^- (\omega )\}} 
 \right] 
\\ 
 & \leq & 
 2 E\left[ 
 \esssup_{\tilde{\omega } \in {\cal N}_\omega} 
 \vert F (\tilde{\omega} ) - F (\omega )\vert^2 
 \right], 
\end{eqnarray*} 
 hence 
$$\frac{1}{\sigma (X)} =\lambda_\infty \geq \frac{(\sqrt{1+h^-_\infty}-1)^2}{2}.$$ 
\end{Proof} 

\begin{prop} 
\label{th9.2} 
 We have 
$$\lambda_\infty 
 = \frac{1}{\sigma (X)} \geq \frac{(\sqrt{h_\infty +1}-1)^2}{4}.$$ 
\end{prop} 
\begin{Proof} 
 We have if $m(F) = 0$: 
\begin{eqnarray*} 
 2 E[\vert D F \vert_\infty ] 
 & \geq & 
 \int_{-\infty}^{+\infty} 
 \pi (\partial \{ F>t\}) dt 
\\ 
& \geq & 
 h_\infty \int_{-\infty}^{+\infty} 
 \min (\pi (\{ F>t\}), 
 \pi ( \{ F\leq t\})) dt 
 = h_\infty E [F]. 
\end{eqnarray*} 
 Applying the above inequality to $(F^+)^2$ we have 
\begin{eqnarray*} 
 h_\infty E [{F^+}^2] & \leq & 
 2 E[ \vert D (F^+)^2 \vert_\infty ] 
\\ 
 & \leq & 
 2 E[\esssup_{\tilde{\omega } \in {\cal N}_\omega}  
  \vert F^+(\omega )- F^+(\tilde{\omega} ) \vert 
 (F^+(\omega ) + F^+(\tilde{\omega} )) ] 
\\ 
 & \leq & 
 2 E[\esssup_{\tilde{\omega } \in {\cal N}_\omega} \vert F^+(\omega )- F^+(\tilde{\omega} ) \vert 
 (F^+(\tilde{\omega} ) - F^+(\omega )) 
\\ 
 & & + 2 \vert F^+(\omega ) - F^+(\tilde{\omega} )\vert F^+ (\omega ) ] 
\\ 
 & \leq & 
 2 E[\esssup_{\tilde{\omega } \in {\cal N}_\omega} 
 (F^+(\omega )- F^+(\tilde{\omega} ) )^2 ] 
\\ 
 & & + 4 E[\esssup_{\tilde{\omega } \in {\cal N}_\omega} 
 \vert F^+(\omega ) - F^+(\tilde{\omega} )\vert F^+ (\omega ) ] 
\\ 
 & \leq & 
 2 E[\esssup_{\tilde{\omega } \in {\cal N}_\omega} 
 (F(\omega )- F(\tilde{\omega} ) )^2 ] 
\\ 
 & & + 4 E[\esssup_{\tilde{\omega } \in {\cal N}_\omega} 
 \vert F(\omega ) - F(\tilde{\omega} )\vert F^+ (\omega ) ]. 
\end{eqnarray*} 
 Similarly we have 
\begin{eqnarray*} 
 h_\infty E [(F^-)^2] & \leq & 
 2 E[\esssup_{\tilde{\omega } \in {\cal N}_\omega} 
 (F(\omega )- F(\tilde{\omega} ) )^2 ] 
 + 4 E[\esssup_{\tilde{\omega } \in {\cal N}_\omega} 
  \vert F(\omega ) - F(\tilde{\omega} )\vert F^- (\omega ) ]. 
\end{eqnarray*} 
 Hence 
\begin{eqnarray*} 
 h_\infty E [F^2] & \leq & 
 h_\infty E [(F^+)^2] + 
 h_\infty E [(F^-)^2] 
\\ 
& \leq & 
 4 E[ \vert D F \vert_\infty^2 ] 
 + 4 E[ \vert D F \vert_\infty \vert F \vert ] 
\\ 
& \leq & 
 4 E[ \vert D F \vert_\infty^2 ] 
 + 
 4 E[ \vert D F \vert_\infty^2 ]^{1/2}
 E[ F^2 ]^{1/2}, 
\end{eqnarray*} 
 which implies 
$$ 
E[ \vert D F \vert_\infty^2 ] \geq 
 E [F^2] 
 \frac{(\sqrt{h_\infty+1}-1)^2}{4}. 
$$ 
 In the general case ($m(F)\not= 0$) we use the fact that $\var F \leq E[(F-m(F))^2]$. 
 Relation \eqref{a5} is proved in Prop.~\ref{9.4}. 
\end{Proof} 
 When $\sigma (X) < \pi/4$, 
 Relation \eqref{a6} also improves the lower bound on 
 $h_\infty$ given in \cite{bobkovgoetze} in the cylindrical 
 (i.e. finite dimensional) case. 
\begin{prop} 
 We have 
\begin{align} 
\label{a1} 
& \lambda_2 = 2 \lambda^+_2 =2\lambda^-_2  = 1, 
\\ 
\label{a2} 
& \lambda_\infty = \frac{1}{\sigma (X)}, 
\\ 
\label{a5} 
& 
 1/\sqrt{2\pi}\leq h_2,  
\\ 
\label{a6} 
& \max \left( 
 \frac{1}{\sqrt{\pi \sigma (X)}} , 
 \frac{1}{2\sigma (X)} \right) 
 \leq h_\infty \leq 
 \frac{4}{\sigma (X)} + \frac{4}{\sqrt{\sigma (X)}}, 
\\ 
\label{a3} 
& \frac{1}{2} \leq h_1 = 2h^+_1 = 2h^-_1 \leq 
 4 + 4 \sqrt{\sigma (X)} 
, 
\\ 
 & \label{a9} 
 h_2^+ 
 \leq 
  \sqrt{1+\sqrt{\sigma (X)}}. 
\\ 
\label{a7} 
& 
 \lambda^+_\infty 
 \leq \frac{h^+_\infty}{2} 
 = \frac{k^+_\infty}{2} 
 \leq \sqrt{\lambda_\infty} = \frac{1}{\sqrt{\sigma (X)}}, 
\\ 
\label{a7.1} 
& 
 \lambda^-_\infty 
 \leq \frac{h^-_\infty}{2} 
 \leq \frac{1}{\sigma (X)} 
 + \sqrt{\frac{2}{\sigma (X)}} 
. 
\end{align} 
\end{prop} 
\begin{Proof} 
\begin{description} 
\item{- Proof of \eqref{a1} and \eqref{a2}.} 
 We have 
\begin{eqnarray*} 
 \var F & \leq & 
 E[\vert DF\vert_{L^2(\sigma)}^2 ] = E[\vert DF\vert_{L^2(\omega )}^2 ] 
 = E\left[\vert DF\vert_{L^2(\frac{\sigma+\omega}{2} )}^2 \right] 
 = 2 E\left[\vert D^\pm F\vert_{L^2(\frac{\sigma+\omega}{2} )}^2 \right] 
\\ 
& \leq & 
 \sigma (X) 
 E[\vert D^\pm F\vert_{L^\infty (\sigma+\omega )}^2 ], 
\end{eqnarray*} 
 hence $\lambda_2 = 2 \lambda^-_2 =2\lambda^+_2  \geq 1$ 
 and $\lambda_\infty \geq 1/\sigma (X)$. 
 Letting $F (\omega ) =\omega (X)$, we have 
$$D_x F = 1_{\{ 
 x\in \omega \}} 
 -1_{\{ 
 x\in \omega^c \}} 
, 
$$ 
 and 
$$ 
 \var (F) = \sigma (X) 
 = E[\vert DF\vert_{L^2(\sigma )}^2] 
 = E[\vert DF\vert_{L^2(\frac{\sigma+\omega}{2})}^2] 
 = \sigma (X) E[\vert DF\vert_{L^\infty (\sigma+\omega )}^2], 
$$ 
 which shows $\lambda_2\leq 1$ 
 and $\lambda_\infty\leq 1/\sigma (X)$. 
\item{- Proof of \eqref{a5}.} From Th.~\ref{9.4}, 
 applying \eqref{e1} to $F=1_A$ we get, 
 since $I (1_A)=0$: 
$$E[ \vert D 1_A\vert_{L^2 (\frac{\sigma+\omega}{2} )} ] 
 \geq 
 \frac{1}{\sqrt{2}} 
 E[ \vert D 1_A\vert_{L^2 (\sigma )} ] 
 \geq 
 \frac{1}{2} 
 I(\pi (A)) \geq \frac{2}{\sqrt{2\pi} } \pi (A) (1-\pi (A)) 
 \geq \frac{1}{\sqrt{2\pi} } \pi (A), 
$$ 
 hence $h_2 \geq 1/\sqrt{2\pi}$. 
\item{- Proof of \eqref{a6}.} 
 We have 
$$E[ \vert D 1_A\vert_{L^\infty (\sigma+\omega )} ] 
 \geq 
 E[ \vert D 1_A\vert_{L^\infty (\sigma )} ] 
 \geq 
 \frac{1}{\sqrt{\sigma (X)}} 
 E[ \vert D 1_A\vert_{L^2 (\sigma )} ] 
 \geq 
 \frac{1}{\sqrt{\pi\sigma (X)} } \pi (A), 
$$ 
 where we used the inequality 
$$I( t ) \geq \sqrt{\frac{2}{\pi}} I_{\rm{var}}(t),$$ 
 with $I_{\rm{var}}(t) = t(1-t)$, $0\leq t \leq 1$, 
 hence $h_\infty \geq 1/\sqrt{\pi \sigma (X)}$. 
 Now if $\pi (A) < 1/2$: 
$$\lambda_\infty \pi (A) \leq 2 \lambda_\infty \pi (A)\pi (A^c ) 
 \leq 2 E[\vert D1_A \vert_{L^\infty (\sigma+\omega )}^2 ] 
 = 2 E[\vert D1_A \vert_{L^\infty (\sigma+\omega )} ], 
$$ 
 hence $\lambda_\infty \leq 2h_\infty$ 
 which, with Prop.~\ref{th9.2} and 
 $h_\infty \geq 1/\sqrt{\pi \sigma (X)}$, 
 proves Relation \eqref{a6}. 
\item{- Proof of \eqref{a3}.} 
 The Clark formula and Lemma~\ref{l1} show that 
 when $\pi (A)\leq 1/2$, 
\begin{eqnarray*} 
\frac{1}{2} \pi (A) & \leq & \var (1_A) \leq E[\vert D1_A\vert_{L^2(\sigma )}^2] 
= E[\vert D1_A\vert_{L^1(\sigma )}] 
\\ 
& = & 2 E\left[\vert D^+1_A\vert_{L^1(\frac{\sigma+\omega}{2} )}\right] 
= 2 E\left[\vert D^-1_A\vert_{L^1(\frac{\sigma+\omega}{2} )}\right] 
= E\left[\vert D 1_A\vert_{L^1(\frac{\sigma+\omega}{2} )}\right] 
\end{eqnarray*} 
 hence 
$$h_1 = 2 h^-_1=2h^+_1 \geq 1/2,$$ 
 which proves the first part of \eqref{a3}. 
 We have 
\begin{eqnarray*} 
 h^+_1 \pi(A) & \leq & 
 E[\vert D^+ 1_A\vert_{L^1 (\frac{\sigma+\omega}{2} )} ] 
 = \frac{1}{2} 
 E[\vert D 1_A\vert_{L^1 (\sigma )} ] 
\\ 
 & \leq &  \frac{1}{2} 
 \sigma (X) E[\vert D 1_A\vert_{L^\infty (\sigma )} ] 
 \leq \frac{1}{2} 
 \sigma (X) E[\vert D 1_A\vert_{L^\infty (\sigma+\omega )} ], 
\end{eqnarray*} 
 hence $h^+_1 \leq \sigma (X) h_\infty/2$, which yields 
 the second part of \eqref{a3} from \eqref{a6}. 
\item{- Proof of \eqref{a9}.} 
 We also have 
\begin{eqnarray*} 
 (h^+_2)^2 \pi(A)^2 & \leq & 
 E\left[\vert D^+ 1_A\vert_{L^2 (\frac{\sigma+\omega}{2} )} \right]^2 
\\ 
 & = &  E\left[1_A \vert D^+ 1_A\vert_{L^2 (\frac{\sigma+\omega}{2} )} \right]^2 
\\ 
 & \leq &  \pi(A) E\left[\vert D^+ 1_A\vert_{L^2 (\frac{\sigma+\omega}{2} )}^2 \right] 
\\ 
 & = &  \pi(A) E\left[\vert D^+ 1_A\vert_{L^1 (\frac{\sigma+\omega}{2} )} \right], 
\end{eqnarray*} 
 hence $(h^+_2)^2 \leq h^+_1$, which proves \eqref{a9}. 
\item{- Proof of \eqref{a7} and \eqref{a7.1}.} 
 Similarly for $\pi (A) \leq 1/2$ we have 
$$\lambda^\pm_\infty \pi (A) \leq 2 \lambda^\pm_\infty \pi (A)\pi (A^c ) 
 \leq 2 E[\vert D^\pm 1_A \vert_{L^\infty (\sigma+\omega )}^2 ] 
 = 2 E[\vert D^\pm 1_A \vert_{L^\infty (\sigma+\omega )} ], 
$$ 
 hence $\lambda^\pm_\infty \leq 2h^\pm_\infty$, 
 and \eqref{a7}, \eqref{a7.1} hold from Prop.~\ref{6.6}. 
\end{description} 
\end{Proof} 
 Clearly the logarithmic Sobolev constants 
$$l^\pm_p = \inf_{0 < \pi (A ) < \frac{1}{2}} 
 \frac{E\left[\vert D^\pm 1_A\vert_{L^p (\frac{\sigma+\omega}{2} )}\right]}{-\pi (A)\log \pi (A) }, 
\quad \mbox{and} 
 \quad 
l_p = \inf_{0 < \pi (A ) < \frac{1}{2}} 
 \frac{ E\left[\vert D 1_A\vert_{L^p (\frac{\sigma+\omega}{2} )}\right]}{-\pi (A)\log \pi (A) } 
$$ 
 vanish, $p\in [1,+\infty ]$, since 
$$ 
l_\infty = \inf_{0 < \pi (A ) < \frac{1}{2}} 
 \frac{\pi (\partial A)}{-\pi (A)\log \pi (A) } 
 = 0, 
$$ 
 (take $A_k = \{\omega (B ) \geq k \}$), i.e. 
 from Th.~\ref{iso} 
 the classical logarithmic Sobolev inequality does not hold on 
 Poisson space. 
 In other terms the optimal constant $\rho_p$ in the inequality 
$$\rho_p \Ent [F^2] \leq E[\vert DF\vert_{L^p (\sigma )}^2]. 
$$ 
 is equal to $0$ for all $p\geq 1$, cf. \cite{ledoux}. 
\section{A remark on Cheeger's inequality} 
\label{s9} 
 This section follows the presentation of \cite{ht} 
 and \cite{stoyanov}, 
 adapting it to the configuration space case. 
 Let $N:\real \rightarrow \real$ be a Young function, i.e. 
 $N$ is convex, even, non-negative, with $N(0)=0$ and $N(x) >0$ 
 for all $x\not= 0$. 
 Let 
$$C_N = \sup_{x>0} \frac{xN'(x)}{N(x)}<\infty. 
$$ 
The Orlicz norm of $F$ is defined as 
$$\Vert F \Vert_N 
 = \inf \left\{ 
 \lambda > 0 \ : \ 
 E\left[ N\left( 
 \frac{F}{\lambda } \right) 
 \right] 
 \leq 1 \right\}. 
$$ 
\begin{theorem}
\label{9.2} 
 For all $F$ such that $m(F)=0$ we have 
$$\Vert F \Vert_N 
 \leq \frac{C_N}{k^+_p} \Vert \vert D F\vert_{L^p (\frac{\sigma+\omega}{2} )}\Vert_N, 
$$ 
 and 
$$E[N(F)] 
 \leq E\left[ 
 N\left( 
 \frac{C_N}{k^+_p} 
 \vert DF\vert_{L^p(\frac{\sigma+\omega}{2} )}\right) 
 \right]. 
$$ 
\end{theorem} 
 For $p=1$ we have $h^+_1=k^+_1$ hence 
$$E[N(F-m(F))] 
 \leq 
 E\left[ 
 N\left( 
 \frac{C_N}{h^+_1} 
 \vert DF\vert_{L^1(\frac{\sigma+\omega}{2} )}\right) 
 \right]. 
$$ 
 If $N(x) = x^p$ we have $C_N = p$ and $\Vert F \Vert_N = \Vert F\Vert_p$, 
 hence for some constant $C(p)$, 
$$C(p ) \Vert F - E[F] \Vert_p 
 \leq \Vert F - m(F) \Vert_p \leq \frac{p}{k_2^+} 
 \Vert \vert DF\vert_{L^2 (\frac{\sigma+\omega}{2} )} \Vert_p. 
$$ 
 For $p=2$ we have $C(2)=1$, hence 
$$\var F \leq \frac{4}{(k_2^+)^2} E\left[\vert DF\vert^2_{L^2(\frac{\sigma+\omega}{2} )}\right], 
$$ 
 and 
$$k_2^+\leq 2.$$ 
 In the particular case $N(x)=x^p$ we have the following better result. 
\begin{theorem}
\label{g1} 
 For all $F$ such that $m(F)=0$ we have 
$$E[\vert F\vert^p] 
 \leq E\left[ 
 \left( 
 \frac{p}{h_1} \vert DF\vert_{L^1(\frac{\sigma+\omega}{2} )} 
 \right)^p  
 \right], 
$$ 
 and 
$$\Vert F \Vert_p 
 \leq \frac{p}{h_1} \Vert \vert D F\vert_{L^1 (\frac{\sigma+\omega}{2} )}\Vert_p. 
$$ 
\end{theorem} 
 We also have the following. 
\begin{prop} 
\label{g2} 
 Let $I_{\rm{var}}(t) = t(1-t)$, $0\leq t \leq 1$ 
 and let $\tilde{b}_p$ denote the optimal constant in the 
 inequality 
$$ 
I_{\rm{var}} 
 (E[F]) \leq E\left[\sqrt{I_{\rm{var}} ( F )^2 + \frac{1}{\tilde{b}_p} 
 \vert DF\vert_{L^2 (\frac{\sigma +\omega}{2} )}^2}\right]. 
$$ 
 We have $\tilde{b}_p \geq 
 \left( 
 1 - \frac{1}{\sqrt{2}} 
 \right) 
 k^+_p. 
$ 
\end{prop} 
\section{Appendix} 
 In this appendix we state the proofs 
 of Th.~\ref{9.2}, Th.~\ref{g1} and 
 Prop.~\ref{g2}, which are based on classical 
 arguments, cf. \cite{bobkovhoudre}, \cite{stoyanov}. 
\begin{Proofx} {\em of Th.~\ref{9.2}}. 
 By the mean value theorem we have 
$$E[\vert D^+ N(F) \vert_{L^p(\frac{\sigma+\omega}{2} )} ] 
 \leq E[N'(F) \vert D^+ F \vert_{L^p(\frac{\sigma+\omega}{2} )} ]. 
$$ 
 On the other hand, if $\Vert F\Vert_N=1$, 
\begin{eqnarray*} 
k^+_p E[N(F)] 
 & = & 
 k^+_p E[N(F^+)] 
 + 
 k^+_p E[N(F^-)] 
\\ 
& \leq & 
E\left[\vert D^+ N(F^+ )\vert_{L^p (\frac{\sigma+\omega}{2} )}\right] 
+ 
E\left[\vert D^+ N(F^- )\vert_{L^p (\frac{\sigma+\omega}{2} )}\right] 
\\ 
& \leq & 
E\left[N'(F^+ ) 
 \vert D^+ F^+\vert_{L^p (\frac{\sigma+\omega}{2} )}\right] 
+ 
E\left[N'(F^- ) 
 \vert D^+ F^- \vert_{L^p (\frac{\sigma+\omega}{2} )}\right] 
\\ 
& \leq & 
E\left[N'(\vert F \vert ) 
 \vert D^+ F \vert_{L^p (\frac{\sigma+\omega}{2} )}\right] 
\\ 
& \leq & 
C_N \Vert \vert D^+ F \vert_{L^p (\frac{\sigma+\omega}{2} )} 
 \Vert_N E[N(F)], 
\end{eqnarray*} 
 where we used the generalization of the H\"older inequality 
$$E[N'(\vert F \vert ) G ] \leq E[N'(\vert F\vert )\vert F\vert]$$ 
 which holds since 
 $1=E[N(G)] \leq E[N(\vert F\vert)]=1$, cf. Lemma~2.1 of \cite{bobkovhoudre}, 
 applied to $\vert F\vert$ and 
$$G= \vert D^+ F \vert_{L^p (\frac{\sigma+\omega}{2} )} 
 \left( 
 \Vert \vert D^+ F \vert_{L^p (\frac{\sigma+\omega}{2} )} 
 \Vert_N 
 \right)^{-1}.$$ 
 Hence 
$$ 
k^+_p  \leq C_N \Vert \vert D^+ F \vert_{L^p (\frac{\sigma+\omega}{2} )} 
 \Vert_N. 
$$ 
 Since $\Vert F \Vert_N =1$, 
 we have 
$$ 
 \Vert F \Vert_N 
 \leq 
 \frac{C_N}{k^+_p} 
 \Vert \vert D^+ F \vert_{L^p (\frac{\sigma+\omega}{2} )} 
 \Vert_N, 
$$ 
 for all $F$ with $m(F)=0$. 
 The second statement is proved by application 
 of the preceding to $N_\alpha (x) = N(x) / \alpha$, $\alpha >0$, 
 as in Th.~3.1 of \cite{bobkovhoudre}. 
\end{Proofx} 
\begin{Proofx} {\em of Th.~\ref{g1}}. 
 We note that 
$$E\left[\vert D\vert F\vert^p \vert_{L^1(\frac{\sigma+\omega}{2} )} \right] 
 \leq p E\left[\vert F\vert^{p-1} \vert D F \vert_{L^1(\frac{\sigma+\omega}{2} )} \right], 
$$ 
 and apply an argument similar to the proof of Th.~\ref{9.2}, 
 with $C_N=p$. 
\end{Proofx} 
\begin{Proofx} {\em of Prop.~\ref{g2}}. 
 The proof is identical to Theorem~4.11 in \cite{stoyanov}. 
 The generalization of Cheeger's inequality applied to 
 $N(x) = \sqrt{1+x^2}-1$ gives $C_N=2$ and 
$$E[N(F)] \leq E\left[ 
 N\left( 
 \frac{2}{k^+_p}\vert DF\vert_{L^p(\sigma+\omega )}\right) 
 \right]. 
$$ 
 We have with 
 $c=\sqrt{2}-1$ and $c_1=k^+_p/2$: 
\begin{eqnarray*} 
 c I_{\rm{var}} (E[F]) 
 & = & 
 c\var (F) + cE[F(1-F)] 
\\ 
 & \leq & c E[F(1-F)] + E[\sqrt{1+F^2}-1] 
\\ 
 & \leq & c E\left[\sqrt{c^2(F(1-F))^2 + \vert DF\vert_{L^p(\frac{\sigma + \omega}{2})}^2/c_1^2}\right], 
\end{eqnarray*} 
 hence 
$$I_{\rm{var}} (E[F]) 
 \leq 
 E[\sqrt{I_{\rm{var}}(F)^2 + \vert DF\vert_{L^p(\frac{\sigma + \omega}{2})}^2/(cc_1)^2}].  
$$ 
\end{Proofx} 

\baselineskip0.5cm

\small 

\def\cprime{$'$} \def\polhk#1{\setbox0=\hbox{#1}{\ooalign{\hidewidth
  \lower1.5ex\hbox{`}\hidewidth\crcr\unhbox0}}}
  \def\polhk#1{\setbox0=\hbox{#1}{\ooalign{\hidewidth
  \lower1.5ex\hbox{`}\hidewidth\crcr\unhbox0}}} \def\cprime{$'$}

\bigskip 

\noindent\sc Laboratoire d'Analyse et de Math\'ematiques Appliqu\'ees, 
CNRS UMR 8050, Universit\'e Paris XII, 94010 Cr\'eteil Cedex, France, and 
School of Mathematics,
Georgia Institute of Technology,
Atlanta, Ga 30332 USA. 
\\ 
{\tt houdre@math.gatech.edu}

\bigskip
\noindent
D\'epartement de Math\'ematiques,
 Universit\'e de La Rochelle, Avenue Michel Cr\'epeau, 
 17042 La Rochelle, France. 
\\ 
{\tt nprivaul@univ-lr.fr}

\end{document}